\documentclass[11pt]{amsart}
\usepackage{amsfonts,amssymb,amsthm,amsmath,amsxtra,amscd,verbatim,eucal}
\usepackage[all]{xy}
\usepackage[dvips]{graphics}

\theoremstyle{plain}
\newtheorem{thm}{Theorem}[section]
\newtheorem{lem}[thm]{Lemma}
\newtheorem{prop}[thm]{Proposition}
\newtheorem{cor}[thm]{Corollary}

\newtheorem{theoremalpha}{Theorem}

\theoremstyle{definition}
\newtheorem{defi}[thm]{Definition}

\theoremstyle{remark}
\newtheorem{rmk}[thm]{Remark}

\newtheorem*{ack}{Acknowledgements}

\def\Z{{\mathbb Z}}

\def\R{{\mathbb R}}
\def\Q{{\mathbb Q}}
\def\P{{\mathbb P}}

\def\E{\mathcal{E}}
\def\FF{\mathcal{F}}

\def\H{\mathcal{H}}

\def\O{\mathcal{O}}

\def\mm{\mathfrak{m}}

\def\a{\alpha}
\def\b{\beta}
\def\g{\gamma}
\def\d{\delta}

\def\f{\phi}
\def\ff{\psi}

\def\l{\lambda}

\def\p{\pi}
\def\r{\rho}
\def\s{\sigma}
\def\t{\tau}

\def\D{\Delta}
\def\G{\Gamma}
\def\LL{\Lambda}
\def\S{\Sigma}

\def\T{\Theta}

\def\.{\cdot}
\def\^{\widehat}
\def\~{\widetilde}
\def\o{\circ}
\def\ov{\overline}
\def\rat{\dashrightarrow}
\def\surj{\twoheadrightarrow}
\def\inj{\hookrightarrow}

\def\({\left(}
\def\){\right)}
\def\ratquot{/ \hskip-2.5pt /}

\DeclareMathOperator{\codim} {codim}

\DeclareMathOperator{\rk} {rk}

\DeclareMathOperator{\Pic} {Pic}

\DeclareMathOperator{\im} {Im}
\DeclareMathOperator{\Hom} {Hom}

\DeclareMathOperator{\Locus} {Locus}
\DeclareMathOperator{\RC} {RC}
\DeclareMathOperator{\Chow} {Chow}
\DeclareMathOperator{\bir} {bir}

\DeclareMathOperator{\CNE} {\ov{NE}}
\DeclareMathOperator{\Nef} {Nef}

\begin{document}

\title{Ample subvarieties and rationally connected fibrations}

\author[M.C. Beltrametti]{Mauro C. Beltrametti}
\address{Dipartimento di Matematica, Universit\`a di Genova,
Via Dodecaneso 35, I-16146 Genova, Italy}
\email{{\tt beltrame@dima.unige.it}}

\author[T. de Fernex]{Tommaso de Fernex}
\address{Department of Mathematics, University of Utah,
155 South 1400 East, Salt Lake City, UT 84112, USA}
\email{{\tt defernex@math.utah.edu}}

\author[A. Lanteri]{Antonio Lanteri}
\address{Dipartimento di Matematica ``F. Enriques'', Universit\`a di Milano,
Via C. Saldini 50, I-20133 Milano, Italy}
\email{{\tt lanteri@mat.unimi.it}}

\thanks{
All authors acknowledge support by MIUR National Research Project
``Geometry on Algebraic Varieties" (Cofin 2004). The research of
the second author was partially supported by NSF grants DMS
0111298 and DMS 0548325. The third author acknowledges partial
support by the University of Milan (FIRST 2003).}




\begin{abstract}
Under some positivity assumptions,
extension properties of rationally connected
fibrations from a submanifold to its ambient variety are studied.
Given a family of rational curves on a complex projective
manifold $X$ inducing a covering family
on a submanifold $Y$ with ample normal bundle in $X$,
the main results relate, under suitable conditions, the associated
rational connected fiber structures on $X$ and on $Y$.
Applications of these results include an extension theorem for
Mori contractions of fiber type and a classification
theorem in the case $Y$ has a structure of projective
bundle or quadric fibration.
\end{abstract}

\maketitle

\section*{Introduction}

Both from a biregular and a birational standpoint,
the geometry of algebraic varieties is often studied in terms of
their fibrations.
Given a smooth complex projective variety $X$,
there are two parallel theories able to detect on $X$
fibrations whose fibers are varieties of negative Kodaira dimension.
The first one, initiated by Mori, associates
a morphism $X \to S$, whose general positive
dimensional fiber is a Fano variety,
to any $K_X$-negative extremal ray of the cone of effective curves of $X$.
The second theory, introduced in works of Campana and of
Koll\'ar, Miyaoka and Mori, produces a rational map $X \rat S$,
with proper and rationally chain connected fibers within its domain
of definition, from any family of rational curves on $X$, the most studied
case being the one in which the family of
all rational curves of $X$ is considered.

The purpose of this paper is to study extension properties of
fibrations of the types described above, when
these fibrations are defined on a smooth subvariety $Y$ of $X$
with ample normal bundle:
our goal is to determine additional conditions that guarantee that
given fiber structures on $Y$ extend to analogous structures on $X$,
and to compare the corresponding fibrations.
Following the terminology of \cite{Hartshorne}, we will refer to such a
$Y$ as an ``ample subvariety'' of $X$.
We remark that, in the special case of codimension one,
this setting is more general than the classical setting
of ample divisors.

Our first result is an extension property for rationally connected
fibrations. The precise formulation requires some additional notation, and
is given in Theorems~\ref{thm:main-RC-extension:part-I}
and~\ref{thm:main-RC-extension:part-II}.
Roughly speaking, Theorem~\ref{thm:main-RC-extension:part-I}
says that if $Y$ is a submanifold
with ample normal bundle,
the inclusion inducing a surjection $N^1(X) \to N^1(Y)$,
and $V$ is a family of rational curves
on $X$ inducing a covering family $V_Y$ on $Y$,
then there is a commutative diagram
$$
\xymatrix{
Y \ar@{-->}[d]_{\p} \ar@{^{(}->}[r] & X \ar@{-->}[d]^{\f} \\
Y\ratquot_{V_Y} \ar[r]^\d & X\ratquot_V,
}
$$
where $\f$ is the rational map associated to the family $V$,
$\p$ is the map associated to $V_Y$, and $\d$ is
a surjective morphism of normal varieties;
here $X\ratquot_V$ and $Y\ratquot_{V_Y}$ denote the
respective rational quotients.

With the second theorem, we determine sufficient conditions
to ensure that the morphism $\d$ in the above diagram is generically finite;
this is one of the core results of the paper
(see Theorem~\ref{thm:main-RC-extension:part-II}).

\begin{theoremalpha}\label{thm:intro:delta-gen-finite}
Assuming that the normal bundle of Y in X is ample and that the inclusion of
$Y$ induces a surjection $N^1(X) \surj N^1(Y)$,
the morphism $\d$ is generically finite if one of the following two
conditions is satisfied:
\begin{enumerate}
\item[(a)]
$V$ is an unsplit family; or
\item[(b)]
$\codim_XY < \dim Y - \dim Y\ratquot_{V_Y}$.
\end{enumerate}
\end{theoremalpha}

When the rational quotient $Y\ratquot_{V_Y}$ is one-dimensional,
it turns out from a general fact that $\p$ is a morphism
(see Proposition~\ref{prop:proper-maps-over-curves-are-morphisms});
using this fact, the previous result can be
improved in this case, and we obtain a commutative diagram as the one above
where now the vertical arrows are morphisms
and $\d$ is a finite morphism between smooth curves
(see Corollary~\ref{cor:RC-extension-unsplit-over-curve}).
The situation when $Y$ is the zero scheme of a regular section of
an ample vector bundle on $X$ and $\p$ is the MRC-fibrations over a
base of positive geometric genus was also treated in
\cite{Occ} (see Remark~\ref{rmk:Occhetta-1}).

Next, we address the problem of extending extremal Mori contractions
of fiber-type. Using the previous result,
we prove the following theorem (see Theorem~\ref{thm:extension-Mori-contraction}).

\begin{theoremalpha}\label{thm:intro:extension-Mori-contr}
Let $Y$ be a submanifold of a projective manifold $X$ with ample normal bundle,
and assume that the inclusion induces an isomorphism $N^1(X) \cong N^1(Y)$.
Let $\p\colon Y \to Z$ be an extremal Mori contraction of fiber-type,
and let $W$ be an irreducible covering family of rational curves on $Y$
that are contracted by $\p$. If these curves do not break in $X$
under deformation (that is, if the family they generate in $X$ is
unsplit), then $\p$
extends to an extremal Mori contraction $\f\colon X \to S$, and
there is a commutative diagram
$$
\xymatrix{
Y \ar[d]_{\p} \ar@{^{(}->}[r] & X \ar[d]^{\f} \\
Z \ar[r]^\d & S,
}
$$
where $\d$ is a finite surjective morphism.
\end{theoremalpha}

In the special case when $Y$ is defined
by a regular section of an ample vector bundle on $X$,
related results were obtained in \cite{AO,dFL,Occ}
(see Remark~\ref{rmk:Occhetta-2}).

In the last section of the paper, the above results are finally applied
to classify, under suitable conditions, projective manifolds
containing ``ample subvarieties'' that admit a structure of
projective bundles or quadric fibrations.
Beginning with classical results on hyperplane sections,
there has been given evidence to the fact that projective manifolds are,
so to speak, at least as ``special'' as their ample divisors
(cf. \cite{Sommese}). The study of projective bundles and quadric fibrations
embedded in projective manifolds as zero schemes of
regular sections of ample vector bundles
was already undertaken in \cite{LM1,LM2,LM3,dF2,Occ},
where classification results were obtained when
the base of the fibration, if positive dimensional,
has positive geometric genus.
Under the additional assumption of a polarization on the ambient variety
inducing a relatively linear polarization on the fibration,
classification results were obtained in \cite{dF1,AO,LM4}.
With some dimensional restrictions
on the fiber structure of the subvariety in terms of its
codimension, we extend such results in two ways.

In our first generalization, we consider fibrations
over curves and drop the hypothesis that
the subvariety is defined by a regular section of an ample vector bundle,
only assuming that the normal bundle is ample.
In this case, since the positivity condition is local near $Y$,
we additionally require that the inclusion induces an isomorphism on the
Picard groups. In the second generalization,
we allow the base of the projective fibration to
have arbitrary dimension.
In either case, we do not put any additional restriction on the base
of the fibration and do not require {\it a priori} any global polarization
(such a polarization will turn out to exist {\it a posteriori}).
We state here the extension property that
we obtain in the case of projective bundles.

\begin{theoremalpha}\label{thm:intro:extension-P-bundle}
Let $Y$ be a submanifold of a projective manifold $X$ with ample normal bundle,
of codimension
\begin{equation}\label{eq:hp:small-codim}
\codim_XY < \dim Y - \dim Z.
\end{equation}
Assume that $Y$ admits a projective bundle
structure $\p \colon Y \to Z$ over a smooth projective
variety $Z$, and that one of the following two situations occurs:
\begin{enumerate}
\item $Z$ is a curve and the inclusion of $Y$ in $X$ induces an
isomorphism $\Pic(X) \cong \Pic(Y)$; or \item $Y$ is the zero
scheme of a regular section of an ample vector bundle on $X$.
\end{enumerate}
Then $\p$ extends to a projective bundle structure $\~\p\colon X \to Z$ on $X$,
giving a commutative diagram
$$
\xymatrix{
Y \ar[d]_\p \ar@{^{(}->}[r] & X \ar[dl]^{\~\p} \\
Z,
}
$$
and the fibers of $\p$ are linearly embedded in the fibers of $\~\p$.
\end{theoremalpha}

A similar statement is also proven in the case of quadric fibrations
with irreducible reduced fibers and relative Picard number one.
For the full result, see Theorem~\ref{thm:P-bundle-Q-fibration}.

The proofs of the above results rely on various properties concerning
families of rational curves on projective manifolds
and their associated rational quotients, which are collected
in Section~2. General facts from deformation theory of rational curves
needed in the sequel are recalled in this section,
and several new notions are introduced and studied.
Among other things, we study numerical properties of
families of rational curves, introducing in particular
the notion of numerically covering family
(see Definition~\ref{def:num-covering-family}),
which turns out to be very useful in our context.
Appropriate notions for the extension and restriction of families
of rational curves are also introduced in relation
to given submanifolds of the ambient variety.

\begin{ack}
We would like to thank C.~Araujo, L.~B\u adescu, C.~Casagrande and P.~Ionescu
for useful discussions, and the referee for many precious
remarks and suggestions.
We are grateful to the Istituto Nazionale di Alta Matematica
and to the University of Milan (FIRST 2003)
for making this collaboration possible.
\end{ack}

\section{Conventions and basic notation}

We work over the complex numbers, and use standard notation in
algebraic geometry, although tensor product between line bundles
is often denoted additively. For $n \ge 1$, we denote by $\Q^n$
the smooth quadric hypersurface of $\P^{n+1}$. We use the word
{\it general} to mean that the choice is made outside a properly
contained (Zariski) closed subset, and the word {\it very general}
to mean that the choice is made outside a countable union of such
subsets. For any projective algebraic variety $X$, we denote by
$N^1(X) := (\Pic(X)/\equiv)\otimes\R$ the N\'eron--Severi space of
$X$, and by $N_1(X) := (Z_1(X)/\equiv)\otimes\R$ the space
generated by numerical classes of curves on $X$. The dimension of
these spaces, which is equal to the Picard number of $X$, is
denoted by $\r(X)$. For a morphism of algebraic varieties $X \to
S$, we set $N_1(X/S) := (Z_1(X/S)/\equiv)\otimes\R$, and we denote
by $\r(X/S)$ the dimension of this space. The numerical class of a
line bundle $L$ (resp., of a curve $C$) on $X$ is denoted by $[L]$
(resp., $[C]$). We denote by $\CNE(X)$ the closure of the cone in
$N_1(X)$ spanned by numerical classes of effective curves, and by
$\Nef(X)$ its dual cone, namely, the cone in $N^1(X)$ spanned by
classes of nef divisors. In the case $X$ is a projective variety
with $\Pic(X) \cong \Z$, we will denote by $\O_X(1)$ the ample
generator of this group.

\section{Families of rational curves on varieties}

Throughout this section we consider morphisms from the projective
line $\P^1$ to a smooth projective variety $X$; for convenience,
we fix once for all two distinct points $0$ and $\infty$ in
$\P^1$. We will use notation as in \cite{KoBook}, some of which we
recall below. For general facts on the theory of
deformation of rational curves on varieties, we refer to
\cite{KoBook,Debarre}.

Let $\Hom(\P^1,X)$ be the scheme parameterizing morphisms from $\P^1$ to $X$.
We will denote a morphism $f \colon \P^1 \to X$ by $[f]$
whenever we think of it as a point of $\Hom(\P^1,X)$.
If $\S \subseteq \P^1$ and $\T \subseteq X$ are closed subschemes,
then we denote by
$$
\Hom(\P^1,X;\S \to \T)
$$
the subscheme of $\Hom(\P^1,X)$
parameterizing those morphisms $f$ such that $f(\S) \subseteq \T$
(the image $f(\S)$ of $\S$ is here intended scheme-theoretically).
In the particular case $\S = \{0\}$ and $\T = \{x\}$ for some
$x \in X$, we also use the notation $\Hom(\P^1,X;0 \mapsto x)$.
Note that, for any closed subscheme $j \colon \S \inj \P^1$,
there is a natural morphism
$$
j^* \colon \Hom(\P^1,X) \to \Hom(\S,X)
$$
defined by $j^*([f]) := [f\o j]$.

\begin{defi}
A morphism $f \colon \P^1 \to X$ is {\it free} if it is non-constant
and $f^*T_X$ is nef.
\end{defi}

We will use the following well-known properties
(see \cite[Theorem~II.1.7]{KoBook}).

\begin{prop}\label{prop:free=smooth}
Let $f \colon \P^1 \to X$ be a non-constant morphism, and let $x = f(0)$.
Let $\mm_0 \subset \O_{\P^1}$ be the maximal ideal sheaf of $0 \in \P^1$.
\begin{enumerate}
\item
The Zariski tangent space of $\Hom(\P^1,X)$ at $[f]$ is
isomorphic to $H^0(\P^1,f^*T_X)$, and $\Hom(\P^1,X)$ is smooth
at $[f]$ if $h^1(\P^1,f^*T_X) = 0$.
\item
Similarly, the Zariski tangent space of
$\Hom(\P^1,X;0 \mapsto x)$ at $[f]$ is
isomorphic to $H^0(\P^1,f^*T_X \otimes \mm_0)$,
and $\Hom(\P^1,X; 0 \mapsto x)$ is smooth
at $[f]$ if $h^1(\P^1,f^*T_X \otimes \mm_0) = 0$.
\end{enumerate}
In particular, if $f$ is free, then
$\Hom(\P^1,X)$ and $\Hom(\P^1,X; 0 \mapsto x)$
are smooth at $[f]$.
\end{prop}

For closed subschemes $\T \subseteq X$ and $\S \subseteq \P^1$,
and any subset $S \subseteq \Hom(\P^1,X)$, we denote
$$
S(\S \to \T) := S \cap \Hom(\P^1,X;\S \to \T),
$$
and, in a similar way, we define $S(0 \mapsto x)$ for any point $x \in X$.

\begin{defi}
Given a subset $S \subseteq \Hom(\P^1,X)$, the image of the
universal map $\P^1 \times S \to X$ is called {\it locus} of $S$,
and is denoted by $\Locus(S)$. Analogous definitions are given for
the loci of $S(\S \to \T)$ and $S(0 \mapsto x)$, which are
respectively denoted by $\Locus(S;\S \to \T)$ and $\Locus(S;0
\mapsto x)$.
\end{defi}

We now restrict our attention to the open subscheme
$\Hom_{\bir}(\P^1,X)$ of $\Hom(\P^1,X)$ parameterizing morphisms
that are birational to their images. The previous notation is
adapted to $\Hom_{\bir}(\P^1,X)$ in the obvious way.

\begin{defi}
Any union $V = \cup_{\a \in A} V_\a$ of irreducible components
$V_\a$ of $\Hom_{\bir}(\P^1,X)$ is said to be a {\it family} of
(parameterized) rational curves on $X$. If $V$ consists of only
one irreducible component of $\Hom_{\bir}(\P^1,X)$, then $V$ is
said to be an {\it irreducible family}. If $\Locus(V_\a)$ is dense
in $X$ for every $\a \in A$, then $V$ is said to be a {\it
covering family}.
\end{defi}

\begin{rmk}
The scheme $\Hom_{\bir}(\P^1,X)$ has at most countably many irreducible
components. In particular, the same holds for
every family of rational curves on $X$.
\end{rmk}

\begin{prop}\label{prop:dense=free}
Let $V$ be an irreducible family of rational curves on $X$. Then the locus
$\Locus(V)$ of $V$ is dense in $X$ if and only if $f$ is free for a general
$[f] \in V$.
\end{prop}

\begin{proof}
$\Locus(V)$ is dense in $X$ if and only if the differential
of the universal map $\P^1 \times V \to X$ has rank equal
to $\dim X$ at a general point $(p,[f]) \in \P^1 \times V$.
By \cite[Proposition~II.3.10]{KoBook}, this occurs if and only
if $f$ is free for a general $[f] \in V$.
\end{proof}

\begin{defi}
Given any closed subset $S \subseteq \Hom_{\bir}(\P^1,X)$,
we denote by $\langle S \rangle$ the union of all irreducible components
of $\Hom_{\bir}(\P^1,X)$ that contain at least one irreducible component
of $S$, and by $\;\rangle S \langle\;$ the union of all irreducible components
of $\Hom_{\bir}(\P^1,X)$ that are irreducible components of $S$.
We call $\langle S \rangle$ the {\it minimal family generated by $S$},
and $\;\rangle S \langle\;$ the {\it maximal family contained in $S$}.
\end{defi}

Analogous notation and definitions will be adopted for $S
\subseteq \Hom(\P^1,X)$.

\begin{defi}
Let $V$ be a family of rational curves on $X$.
We say that a curve $C \subset X$ is a {\it $V$-curve}
if $C = f(\P^1)$ for some $[f] \in V$.
A {\it $V$-chain of length $\ell$} is the union of $\ell$ distinct curves
$f_i(\P^1) \subset X$ ($1 \le i \le \ell$)
parameterized by elements $[f_i] \in V$ such that
$f_{i+1}(0) = f_i(\infty)$ for every $1 \le i \le \ell-1$.
\end{defi}

Consider a family $V = \cup V_\a$ of rational curves on $X$.
Let $V_\a'$ be the irreducible components of the
normalization $\Hom_{\bir}^n(\P^1,X)$ of $\Hom_{\bir}(\P^1,X)$ that map to $V_\a$,
and let $\ov{V}$ be the closure in $\Chow_1(X)$
of the image of $\cup V_\a'$ via the
natural morphism
\begin{equation}\label{eq:Hom-to-Chow}
\Hom_{\bir}^n(\P^1,X) \to \Chow_1(X)
\end{equation}
(see \cite[Comment~II.2.7]{KoBook}).
By \cite[Theorem~I.5.10]{KoBook}, $\ov{V}$ is proper. Associated to $\ov{V}$
(hence to $V$),
one can define a {\it proper proalgebraic relation} on $X$, as explained in
\cite[Section~IV.4]{KoBook} (see in particular \cite[Example~IV.4.10]{KoBook}).
We denote this relation by $\RC_V$.
We have that $(x,y) \in \RC_V$ for two very general
points $x, y \in X$ if and only
if there is a $V$-chain containing $x$ and $y$.

\begin{defi}
$X$ is said to be {\it $\RC_V$-connected} if there is only one $\RC_V$-equivalence
class. We say that $X$ is {\it $V$-chain connected} if every two
points of $X$ lie on a $V$-chain.
More generally, a closed subset $T\subset X$ is said to
be {\it $V$-chain connected}
if any two points of $T$ lie on a $V$-chain supported on $T$.
\end{defi}

Clearly if $X$ is $V$-chain connected, then it is $\RC_V$-connected.

\begin{defi}
A family of rational curves $V$ is said to be
an {\it unsplit family} if it is irreducible and, after normalization, its image
in $\Chow_1(X)$ via the map \eqref{eq:Hom-to-Chow} is a proper scheme.
\end{defi}

\begin{rmk}\label{rmk:HC=1-implies-V-unsplit}
If $V$ is an irreducible family of rational curves on $X$, and
there is an ample vector bundle $\E$ on $X$ such that $\det\E\.C <
2\rk \E$ for a $V$-curve $C$, then $V$ is unsplit. Indeed, if
$\sum_{i=1}^k [C_i]$ is any degeneration of the cycle $[C]$, then
$\det\E\.C = \sum_{i=1}^k \det\E\.C_i \ge k \rk\E$ by the
ampleness of $\E$. Thus $k = 1$ and $[C_1]$ lies in the image of
the normalization of $V$ in $\Chow_1(X)$ via the map
\eqref{eq:Hom-to-Chow}.
\end{rmk}

\begin{prop}\label{prop:unsplit-rho=1}
If $X$ is $\RC_V$-connected by a family $V = \cup V_\a$
with each $V_\a$ unsplit, then
it is $V$-chain connected and $N_1(X)$ is generated by
numerical classes of $V$-curves. In particular, if $V$ is unsplit, then
$X$ has Picard number $\r(X) = 1$.
\end{prop}

\begin{proof}
Fix $\ell \gg 0$ such that every two very general points of $X$ are
connected by a $V$-chain of length $\ell$,
and let $S \subset \Chow_1(X)$ be the subset parameterizing
connected 1-cycles with $\ell$ components (counting multiplicities)
supported on $V$-curves. Since each component $V_\a$
of $V$ is unsplit, it follows that $S$ is proper.
By the choice of $\ell$, there is an irreducible component $T$ of $S$ such that
$u^{(2)}\colon U\times_T U \to X \times X$ is dominant,
where $U \to T$ is the universal family and $u \colon U \to X$
is the cycle map. Note that $T$ is
a closed subvariety of $\Chow_1(X)$, and in particular it is proper.
Therefore the image of $u^{(2)}$ is proper,
and thus $u^{(2)}$ is surjective. This implies that $X$
is $V$-chain connected. The second assertion follows from
\cite[Proposition~IV.3.13.3]{KoBook}.
\end{proof}

If $V = \Hom_{\bir}(\P^1,X)$, then it is a result of Campana
\cite{Campana} and Koll\'ar--Miyaoka--Mori \cite{KMM} that one can
``pass to the quotient''. The natural generalization of this
result when $V$ is an arbitrary family of rational curves was
later given by Koll\'ar (see~\cite[Theorem~IV.4.17]{KoBook}), and
can be stated as follows.

\begin{thm}\label{thm:RC-fibration}
Let $V$ be a family of rational curves on $X$.
Then there is an open set $X^\o \subseteq X$, a
normal projective variety $X\ratquot_V$, and a dominant
morphism $X^\o \to X\ratquot_V$ with connected fibers
and proper over the image,
whose very general fibers are $\RC_V$-equivalence classes in $X$.
\end{thm}

\begin{proof}
We refer the reader to \cite[Theorem~IV.4.17]{KoBook} for the
existence of a proper morphism $X^\o \to Y^\o$ onto a variety $Y^\o$
with connected fibers and
whose very general fibers are $\RC_V$-equivalence classes in $X$.
Since $X^\o$ is normal and the fibers are connected, $Y^\o$ is normal.
It follows by construction that $Y^\o$ is quasi-projective.
Then we can take $X\ratquot_V$ to be the normalization of the closure of $Y^\o$
in some projective embedding.
\end{proof}

\begin{defi}
With the notation as in the previous theorem,
the resulting rational map $X \rat X\ratquot_V$
is called the {\it $\RC_V$-fibration}
of $X$, and $X\ratquot_V$ is the {\it $\RC_V$-quotient}. We also
say that $X\ratquot_V$ (resp., $X \rat X\ratquot_V$) is the rational
quotient (resp., the rational fibration) {\it defined} by $V$. We denote by
$$
\dim\RC_V := \dim X - \dim X\ratquot_V
$$
the dimension of a very general $\RC_V$-equivalence class (which is the
same as the dimension of a general fiber of $X \rat X\ratquot_V$).
\end{defi}

\begin{rmk}
We have $\dim\RC_V > 0$ if and only if $V$ is a covering family.
We would like to stress that $X\ratquot_V$ is
well defined only up to birational equivalence;
throughout the paper, we will often make suitable choices
of rational quotients.
We remark that, after possibly shrinking further $X^\o$, one can always
take a smooth projective model for the rational quotient $X\ratquot_V$.
\end{rmk}

\begin{rmk}
A case of particular interest is when $V = \Hom_{\bir}(\P^1,X)$.
In this case $X\ratquot_V$ is called the {\it {\rm MRC}-quotient}
({\it maximal rational quotient}) of $X$. It follows from a result
of Graber, Harris and Starr \cite{GHS} that in this situation
$X\ratquot_V$ is not uniruled.
\end{rmk}

Every connected subset $S \subseteq \Hom_{\bir}(\P^1,X)$
determines in a natural way a numerical class $[S] \in N_1(X)$, by
taking the class of the curve $f(\P^1)$ for an arbitrary $[f] \in
S$. More generally, we give the following definition.

\begin{defi}
For any subset $S \subseteq \Hom_{\bir}(\P^1,X)$, we denote
$$
\R_{\ge 0}[S] := \sum_{[f] \in S}\R_{\ge 0}[f(\P^1)] \subset N_1(X).
$$
We call $\R_{\ge 0}[S]$ the {\it cone numerically spanned by $S$}.
If $S$ is connected, then we call $[S]$ the {\it numerical class} of $S$.
\end{defi}

\begin{defi}\label{def:num-covering-family}
A family of rational curves $V = \cup_{\a \in A} V_\a$ is said to be
{\it numerically covering} if
there is a subset $A' \subseteq A$ such that
$$
\R_{\ge 0}[V] = \sum_{\a \in A'} \R_{\ge 0}[V_\a]
$$
and $\Locus(V_\a)$ is dense in $X$ for every $\a \in A'$. If $V$
is a numerically covering family, then we denote by $V_{\rm cov}$
the subfamily of $V$ consisting of all the covering irreducible
components of $V$.
\end{defi}

\begin{rmk}
Note that by definition $\R_{\ge 0}[V_{\rm cov}] = \R_{\ge 0}[V]$
for every numerically covering family $V$ of rational curves on
$X$. Clearly, in the notation of the definition, we have $\cup_{\a
\in A'} V_\a \subseteq V_{\rm cov}$.
\end{rmk}

\begin{prop}\label{prop:numerical-condition-on-RC}
Let $V = \cup_{\a \in A} V_\a$ and $W = \cup_{\b \in B} W_\b$ be two
numerically covering families of rational curves on $X$, and suppose that
$$
\R_{\ge 0}[V] \subseteq \R_{\ge 0}[W].
$$
Then the
$\RC_W$-fibration $\ff \colon X \rat X\ratquot_W$ factors through
the $\RC_V$-quotient of $X$, i.e., we have a commutative diagram
$$
\xymatrix{
X\ar@{-->}[d]_{\f} \ar@{-->}[dr]^{\ff} \\
X \ratquot_V \ar@{-->}[r] & X \ratquot_W.
}
$$
In particular, if $\R_{\ge 0}[V] = \R_{\ge 0}[W]$,
then $V$ and $W$ define the same rational fibration.
\end{prop}

\begin{proof}
Consider a resolution of the indeterminacies of $\ff$
$$
\xymatrix{
X \ar@{-->}[dr]_\ff & X' \ar[l]_\s \ar[d]^{\ff'} \\
& X\ratquot_W,
}
$$
and let $E \subset X'$ be the exceptional locus of $\s$.
Note that the locus of indeterminacies of $\ff$
is contained in $X\setminus X^\o$, and hence, by Theorem~\ref{thm:RC-fibration},
it does not dominate $X\ratquot_W$, since it does not meet the very general
$\RC_V$-equivalence class. Therefore
we can assume that $\ff'(E)$ is a proper closed subset of $X\ratquot_W$.

\begin{lem}\label{lem:C-mapped-to-a-point}
If $C \subset X$ is a $V$-curve not contained in $\s(E)$, then its
proper transform $C' \subset X'$ is mapped to a point in
$X\ratquot_W$ by $\ff'$.
\end{lem}

\begin{proof}[Proof of the lemma.]
Let $B' \subset B$ be the subset parameterizing the
irreducible components of $W$ that are covering, so that
$W_{\rm cov} = \cup_{\b \in B'} W_\b$. Since
$$
\R_{\ge 0}[V] \subseteq \R_{\ge 0}[W] = \R_{\ge 0}[W_{\rm cov}],
$$
we can pick general elements $[g_\b] \in W_\b$ for $\b \in B'$, and
numbers $\l_\b \ge 0$, such that, denoting $\G_\b = g_\b(\P^1)$, we have
$$
[C] = \sum_{\b \in B'} \l_\b [\G_\b] \quad\text{in}\quad N_1(X).
$$
By the definition of $\ff$ and Theorem~\ref{thm:RC-fibration} and
the fact that the $W_\b$ are covering families for every $\b \in
B'$, we can assume that $\G_\b \cap \s(E) = \emptyset$ for every
$\b$. Let $C'$ and $\G_\b'$ be the proper transforms of $C$ and
$\G_\b$ on $X'$, and let $D \subset X'$ be the pull-back, via
$\ff'$, of an ample divisor on $X\ratquot_W$. Then $\ff'_*[\G_\b']
= 0$, so we have $\s_*D\.\G_\b = D \.\G_\b' = 0$ for all $\b$, and
thus
$$
D\.C' \le \s_*D \. C = \sum_{\b \in B'} \l_\b (\s_*D\.\G_\b) = 0.
$$
Since $D$ is the pull-back of an ample divisor on $X\ratquot_W$,
this implies that $\ff'_*[C'] = 0$.
\end{proof}

We can now conclude the proof of the proposition. We consider a
very general fiber $G$ of $\f$, and fix two very general points
$x,y \in G$ that are not contained in $\s(E)$. By the generality
of the choices, we can assume that both $x$ and $y$ are outside
the locus of indeterminacies of $\ff$, hence $\ff(x), \ff(y) \not
\in \ff'(E)$. Furthermore, we can assume that $x$ and $y$ are
connected by a chain of $V$-curves. By
Lemma~\ref{lem:C-mapped-to-a-point}, every irreducible component
of this chain that is not contained in $\s(E)$ is mapped to a
point by $\ff$. Therefore, since $\ff(x) \not \in \ff'(E)$, we
deduce that this chain is in fact disjoint from $\s(E)$. Thus we
can apply the lemma to each one of its components. Since the chain
is connected, this implies
$$
\ff(x) = \ff(y).
$$
In view of the generality of the choice of $x$ and $y$ in $G$,
we conclude that there exists a natural rational
map $X\ratquot_V \rat X\ratquot_W$
commuting with the respective projections.
\end{proof}

\begin{defi}
Let $V \subseteq \Hom_{\bir}(\P^1,X)$ be any family of rational curves, and
let $a > 0$ be an integer. We denote by $\|V\|$
the largest family of rational curves on $X$ such that
$\R_{\ge 0}[\|V\|] = \R_{\ge 0}[V]$. We call
$\|V\|$ the family {\it numerically generated} by $V$.
\end{defi}

\begin{rmk}
Quite obviously, $V$ is a subfamily of $\|V\|$, and if $V$ is a
numerically covering family then so is $\|V\|$.
\end{rmk}

The previous proposition implies the following useful property.

\begin{cor}\label{cor:V-versus-dV}
Let $V$ be a numerically covering family of rational curves on $X$.
Then the families $V$, $V_{\rm cov}$, and $\|V\|$
define the same rational fibration.
\end{cor}

\begin{proof}
Since
$\R_{\ge 0}[\|V\|] = \R_{\ge 0}[V] =  \R_{\ge 0}[V_{\rm cov}]$,
it follows from Proposition~\ref{prop:numerical-condition-on-RC}.
\end{proof}

An important theorem, due to Koll\'ar, Miyaoka and Mori, says that
a smooth projective manifold is maximally rationally chain
connected if and only if it is maximally rationally connected (see
\cite[Theorem~IV.3.10]{KoBook}). The proof of this property can be
adapted to the situation in which an arbitrary family of rational
curves $V$ is considered. More precisely, one can prove that if
$X$ is $V$-chain connected, then every two very general points of
$X$ are connected by a $\|V\|$-curve. In the sequel, we will need
the following slightly different version of this property (the
proofs of the two properties are almost the same).

\begin{prop}\label{prop:V-chain-connected=Locus(V,y)-dense}
Let $V$ be a family of rational
curves on $X$, and assume that $X$ is $V$-chain connected.
Let $y \in X$ be any point that is connected by a $V$-chain to a
very general point of $X$. Then
$\Locus(\|V\|,0 \mapsto y)$ is dense in $X$.
\end{prop}

\begin{proof}
Fix a very general point $x \in X$, and let $C = C_1 +\dots+C_n$
be a $V$-chain connecting $x$ to $y$, with $C_i = f_i(\P^1)$.
Let $p_0 = x$, $p_n = y$, and $p_i = f_i(\infty) = f_{i+1}(0)$
for $1 \le i \le n-1$. From here,
we follow step by step the proof of \cite[Complement~IV.3.10.1]{KoBook},
proving inductively that there is a free rational curve $g_i\colon \P^1 \to X$
which connects $p_{i-1}$ to $p_i$;
the resulting chain of free rational curves
is then smoothed into a free rational curve $h \colon \P^1 \to X$
connecting $x$ to $y$. The construction shows that
each $g_{i+1}$ is obtained by smoothing a comb on $f_{i+1}$ with teeth
assembled out of deformations of $g_i$ (see \cite[Definition~II.7.7]{KoBook}).
It follows that $h$ is a $\|V\|$-curve.
\end{proof}

Consider now a submanifold $Y$ of $X$. The inclusion
$i \colon Y \inj X$ naturally induces an injective morphism
$$
i_* \colon \Hom(\P^1,Y) \inj \Hom(\P^1,X),
$$
which is defined by $i_*([g]) := [i\o g]$.
Similar notation will be used for the injection
$\Hom_{\bir}(\P^1,Y) \inj \Hom_{\bir}(\P^1,X)$.
Statements analogous to the following two propositions also
hold for $\Hom_{\bir}(\quad )$ in place of $\Hom(\quad)$.

\begin{prop}\label{prop:free-on-Y=free-on-X}
With the above notation, assume that the normal bundle of $Y$ in
$X$ is ample. Then, for every $[g] \in \Hom(\P^1,Y)$ with $g$
free, the morphism $f := i\o g$ is free on $X$, and the schemes
$\Hom(\P^1,X)$ and $\Hom(\P^1,X;0 \mapsto g(0))$ are smooth at
$[f]$.
\end{prop}

\begin{proof}
The bundle $g^*T_Y$ is nef, since $[g]$ is free.
Thus, by the ampleness of $N_{Y/X}$ and the exact sequence
$$
0 \to g^*T_Y \to f^*T_X \to g^*N_{Y/X} \to 0,
$$
we conclude that $f$ is free.
Then the last two assertions follow from Proposition~\ref{prop:free=smooth}.
\end{proof}

We will need the following generalization of
Proposition~\ref{prop:free=smooth}.

\begin{prop}\label{prop:tangent-sp}
Let $Y$ be a submanifold of a smooth projective variety $X$ and
let $i \colon Y \inj X$ be the inclusion. Let $g \colon \P^1 \to
Y$ be a free morphism, and let $f := i\o g$. Then the Zariski
tangent space of $\Hom(\P^1,X;\{0\} \to Y)$ at $[f]$
sits naturally in an exact sequence
$$
0 \to H^0(\P^1,g^*T_Y) \to T_{[f]}\Hom(\P^1,X;\{0\} \to Y)
\to H^0(\P^1,g^*N_{Y/X} \otimes \mm_0).
$$
Moreover, if the normal bundle of $Y$ in $X$ is ample, then
the sequence completes to a short exact sequence
$$
0 \to H^0(\P^1,g^*T_Y) \to T_{[f]}\Hom(\P^1,X;\{0\} \to Y)
\to H^0(\P^1,g^*N_{Y/X} \otimes \mm_0) \to 0
$$
and $\Hom(\P^1,X;\{0\} \to Y)$ is smooth at $[f]$.
\end{prop}

\begin{proof}
The natural maps $ f^*T_X \to g^*N_{Y/X} \to g^*N_{Y/X}|_{\{0\}}$,
passing to cohomology and taking into account the isomorphism
$T_{[f]}\Hom(\P^1,X) \cong H^0(\P^1,f^*T_X)$ given by
Proposition~\ref{prop:free=smooth}, yield a natural homomorphism
$$
r: T_{[f]}\Hom(\P^1,X) \to (g^*N_{Y/X})|_{\{0\}}.
$$
We have a Cartesian square
$$
\xymatrix{
\Hom(\P^1,X;\{0\} \to Y) \ar[d] \ar[r] & \Hom(\P^1,X) \ar[d] \\
\Hom(\{0\},Y) \ar[r] & \Hom(\{0\},X). }
$$
Note that there are natural identifications $\Hom(\{0\},Y) = Y$ and
$\Hom(\{0\},X) = X$. The Zariski tangent space of $\Hom(\P^1,X;\{0\} \to
Y)$ at $[f]$, viewed as a subspace of $T_{[f]}\Hom(\P^1,X)$, is
contained in ${\rm ker}(r)$. This simply follows from the fact
that the morphism $\Hom(\P^1,X;\{0\} \to Y) \to \Hom(\{0\},X)$
factors through $\Hom(\{0\},Y)$, by taking tangent spaces and
recalling the above natural identifications. Now, the freeness of
$g$ implies that $H^1(\P^1,g^*T_Y) = 0$, hence we get the exact
sequence
$$
0 \to H^0(\P^1,g^*T_Y) \to T_{[f]}\Hom(\P^1,X) \to
H^0(\P^1,g^*N_{Y/X}) \to 0.
$$
This sequence restricts to an exact sequence
\begin{equation}\label{eq:10}
0 \to H^0(\P^1,g^*T_Y) \to {\rm ker}(r)
\to H^0(\P^1,g^*N_{Y/X} \otimes \mm_0) \to 0,
\end{equation}
and the first assertion follows by observing that
$T_{[f]}\Hom(\P^1,X;\{0\}\to Y)$, viewed as a subspace of ${\rm ker}(r)$,
contains the image of $H^0(\P^1,g^*T_Y)$.

Suppose now that $N_{Y/X}$ is ample. By
Proposition~\ref{prop:free-on-Y=free-on-X}, $\Hom(\P^1,X;0 \mapsto
g(0))$ is smooth, hence of dimension $h^0(\P^1,f^*T_X \otimes
\mm_0)$, at $[f]$. Let $V$ be the irreducible component of
$\Hom(\P^1,X)$ containing $[f]$. Note that
$T_{[f]}\Hom(\P^1,X;\{0\} \to Y) = T_{[f]}V(\{0\}\to Y)$. Moreover,
there is a dominant morphism $V(\{0\} \to Y) \to Y$, defined by
$[h] \mapsto h(0)$, whose fiber over a point $y \in Y$ is $V(0
\mapsto y)$. This implies that $\dim V(\{0\} \to Y) = \dim V(0
\mapsto y) + \dim Y$ for a general $y \in Y$, and therefore
\begin{equation}\label{eq:dim2}
\dim V(\{0\} \to Y) = h^0(\P^1,f^*T_X \otimes \mm_0) + \dim Y
\end{equation}
by Proposition~\ref{prop:free=smooth},(b). On the other hand, we
have
\begin{align*}
h^0(\P^1,g^*N_{Y/X} \otimes \mm_0)
&= h^0(\P^1,f^*T_X \otimes \mm_0) - h^0(\P^1,g^*T_Y \otimes \mm_0) \\
&= h^0(\P^1,f^*T_X \otimes \mm_0) - h^0(\P^1,g^*T_Y) +
h^0(g^*T_Y|_{\{0\}}).
\end{align*}
Therefore, by \eqref{eq:10} and $h^0(g^*T_Y|_{\{0\}}) = \dim Y$,
we get
\begin{equation}\label{eq:tangent-sp}
\dim {\rm ker}(r) = h^0(\P^1,f^*T_X \otimes \mm_0) + \dim Y.
\end{equation}
Then, comparing~\eqref{eq:tangent-sp} with~\eqref{eq:dim2}, we
conclude at once that $T_{[f]}\Hom(\P^1,X;\{0\} \to Y)= {\rm ker}(r)$
and that $V(\{0\} \to Y)$ is smooth at $[f]$.
\end{proof}

\begin{defi}\label{def:restriction-of-V}
If $W \subseteq \Hom_{\bir}(\P^1,Y)$ is a family of rational curves on $Y$,
then we call $\langle i_*(W) \rangle \subseteq \Hom_{\bir}(\P^1,X)$
the {\it extension} of $W$ to $X$.
Conversely, for every family of rational curves $V$ on $X$, we call
$\;\rangle i_*^{-1}(V) \langle\; \subseteq \Hom_{\bir}(\P^1,Y)$
the {\it restriction} of $V$ to $Y$.
\end{defi}

\begin{rmk}
If $V$ is an irreducible family on $X$, then
$\;\rangle i_*^{-1}(V) \langle\;$ needs not be irreducible. In fact, the example
of the restriction to a smooth quadric $\Q^2$ of the family of lines in $\P^3$
shows that, in general, different elements in $\;\rangle i_*^{-1}(V) \langle\;$
may even define linearly independent numerical classes in $N_1(Y)$.
In particular, if $V$ is an unsplit family on $X$, then its restriction
to $Y$ may not be an unsplit family.
Although we do not have examples,
it seems likely that $i_*^{-1}(V)$ might fail to be
a family on $Y$ even if $V$ is a family on $X$.
\end{rmk}

We close this section with the following ``relative'' version of
Proposition~\ref{prop:unsplit-rho=1}.
The proof, which can be found within the proof
of \cite[Lemma~1.4.5]{BSW}, uses a ``non-breaking lemma''
due to Wi\'sniewski \cite{Wis}.

\begin{prop}\label{prop:relative-unspli=picard1:occhetta}
Let $V$ be an unsplit family on $X$. Let $Y \subset X$ be a
subvariety, and assume that $\Locus(V,\{0\} \to Y)$ is dense in
$X$. Then, for every curve $\G \subset X$, we have
$$
[\G] = a [\G_Y] + b [C] \quad\text{in}\quad N_1(X),
$$
where $C$ is a $V$-curve, $\G_Y$ is a curve contained in $Y$,
and $a \ge 0$.
\end{prop}

\section{Extension of rationally connected fibrations}

In this section we study
the relationship between rational connected fibrations of a
projective manifold $X$ and those of an ample submanifold $Y \subset X$.
Our main interest is in the case when it is given
on $X$ an irreducible family of rational curves $V$
whose restriction $V_Y$ to $Y$
contains a covering component. Given this situation,
the goal is to compare the associated rationally connected
fibrations $X \rat X\ratquot_V$ and $Y \rat Y\ratquot_{V_Y}$.
However, for technical reasons that will be evident in the
proof of Theorem~\ref{thm:main-RC-extension:part-II},
it is convenient to consider from the beginning a more
general setting, allowing $V$ to be reducible.
The conditions given in items~(i) and~(ii) in the following
two theorems capture the essential properties of $V$ needed
in our arguments.

\begin{thm}\label{thm:main-RC-extension:part-I}
Let $X$ be a smooth projective variety, and assume that $Y \subset
X$ is a smooth subvariety with ample normal bundle. Let $i \colon
Y \inj X$ be the inclusion, and suppose that the
induced map $i^*\colon N^1(X) \to N^1(Y)$ is surjective.
Let $V = \cup_{\a \in A} V_\a$ be a
family of rational curves on $X$, and assume that there is a
subset $B \subseteq A$ such that
\begin{enumerate}
\item[(i)] $\R_{\ge 0}[V] = \sum_{\b \in B} \R_{\ge 0}[V_\b]$, and
\item [(ii)] $\Locus\big(i_*^{-1}(V_\b)\big)$ is dense in $Y$ for
every $\b \in B$.
\end{enumerate}
Let $V_Y := \;\rangle i_*^{-1}(V) \langle\;$
be the restriction of $V$ to $Y$. Then both $V$ and $V_Y$
are numerically covering families (respectively, on $X$ and on $Y$),
and, for suitable choices of the rational quotients,
there is a commutative diagram
$$
\xymatrix{
Y \ar@{-->}[d]_{\p} \ar@{^{(}->}[r] & X \ar@{-->}[d]^{\f} \\
Y\ratquot_{V_Y} \ar[r]^\d & X\ratquot_V,
}
$$
where $\p$ and $\f$ are the projections to the respective
rational quotients and $\d$ is a surjective morphism.
\end{thm}

\begin{rmk}
If one assumes that $V$ is irreducible,
then the hypothesis that $i^* \colon N^1(X) \to N^1(Y)$ be surjective
is unnecessary.
\end{rmk}

\begin{rmk}
An analogous property has been independently observed to hold
in the case $V = \Hom_{\bir}(\P^1,X)$ by A.J.~de~Jong and J.~Starr.
\end{rmk}

The proof of Theorem~\ref{thm:main-RC-extension:part-I}
is based upon the following two lemmas.

\begin{lem}\label{lem:B-dense}
With the same notation and assumptions as in
Theorem~\ref{thm:main-RC-extension:part-I},
$\Locus(V_\b ; \{0\} \to Y)$ is dense in $X$
for every $\b \in B$.
\end{lem}

\begin{proof}
To keep the notation light, we suppose throughout the proof that
$V$ is irreducible, so that $V_\b = V$. By hypothesis,
$\Locus(V_Y)$ is dense in $Y$. Fix a general element $[g]$ of an
irreducible component of $V_Y$ with dense locus in $Y$, and let $f
:= i \o g$. Note that $[f] \in V(\{0\} \to Y)$. Since the chosen
component of $V_Y$ has dense locus in $Y$, we can assume that $y := f(0)$
is a general point of $Y$. We know that $g$ is free by
Proposition~\ref{prop:dense=free}. Thus, since $V$ is an
irreducible component of $\Hom(\P^1,X)$ containing $[f]$,
Proposition~\ref{prop:tangent-sp} applies to say that $V(\{0\} \to
Y)$ is smooth at $[f]$, with Zariski tangent space
$$
T_{[f]} V(\{0\} \to Y) \cong
H^0(\P^1,g^*T_Y) \oplus H^0(\P^1,g^*N_{Y/X} \otimes \mm_0)).
$$
We can view the right hand side as a vector subspace $\LL
\subseteq H^0(\P^1,f^*T_X)$ via the isomorphism \eqref{eq:10}
given in the proof of Proposition~\ref{prop:tangent-sp}. It
follows from the ampleness of $g^*N_{Y/X}$ that $f^*T_X$ is
spanned by the sections in $\LL$ at every point $q \in \P^1
\setminus \{0\}$. By \cite[Proposition~I.2.19]{KoBook}, we
conclude that the differential of the universal map
$$
\P^1 \times V(\{0\} \to Y) \to X
$$
has rank equal to $\dim X$ at $(q,[f])$ for every $q \in \P^1 \setminus \{0\}$.
Since $V(\{0\} \to Y)$ is smooth at $[f]$,
this implies that its locus in $X$ is dense.
\end{proof}

\begin{lem}\label{lem:V-covering-and-V_Y-family}
With the same notation and assumptions as in
Theorem~\ref{thm:main-RC-extension:part-I},
for every $\b \in B$ the family $V_\b$
is covering (in $X$) and every irreducible component
of $i_*^{-1}(V_\b)$ with dense locus in $Y$ is a family (on $Y$).
\end{lem}

\begin{proof}
By the hypothesis~(ii) of
Theorem~\ref{thm:main-RC-extension:part-I}, we can fix an
arbitrary irreducible component $T$ of $i_*^{-1}(V_\b)$ with dense
locus in $Y$. Let $[g] \in T$ be a general element. We know by
Proposition~\ref{prop:dense=free} that $g$ is free. Then by
Proposition~\ref{prop:free-on-Y=free-on-X}, $f := i\o g$ is free
and $\Hom(\P^1,X)$ is smooth at $[f]$, and thus, in particular,
$V_\b$ is a covering family by Proposition~\ref{prop:dense=free}
again. Moreover, we deduce by the smoothness of $\Hom(\P^1,X)$ at
$[f]$ that for every irreducible one-parameter family $[g_t]$ in
$\Hom_{\bir}(\P^1,Y)$ specializing to $[g]$, the corresponding
family $[f_t] := [i \o g_t]$ is contained in $V$. Therefore $[g_t]
\in i_*^{-1}(V)$, and in fact $[g_t] \in T$ by the generality of
the choice of $[g]$ in the component $T$. We conclude that $T$ is
an irreducible component of $\Hom_{\bir}(\P^1,Y)$, and hence, in
particular, of $V_Y$.
\end{proof}

\begin{proof}[Proof of Theorem~\ref{thm:main-RC-extension:part-I}.]
It follows from Lemma~\ref{lem:V-covering-and-V_Y-family} and the hypothesis~(i)
that $V$ is a numerically covering family on $X$ and similarly,
using the fact that the $i_* \colon N_1(Y) \to N_1(X)$ is injective,
it follows that $V_Y$ is a numerically covering family on $Y$.
Let $\f \colon X \rat Z \ratquot_V$ and $\p \colon Y \rat Y\ratquot_{V_Y}$
be the rational quotients, as in Theorem~\ref{thm:RC-fibration}.
Let $X^\o \subseteq X$, $S^\o \subseteq X\ratquot_V$, $Y^\o \subseteq Y$,
and $Z^\o \subseteq Y\ratquot_{V_Y}$ be open subsets such that the maps
$\f$ and $\p$ restrict to proper surjective morphisms
$$
\f^\o \colon X^\o \to S^\o \quad\text{and}\quad
\p^\o \colon Y^\o \to Z^\o.
$$
Let $G$ be a very general fiber of $\f^\o$, and let $x$ be a
general point of $G$. By Lemma~\ref{lem:B-dense}, we can assume
that $x \in \Locus(V;\{0\} \to Y)$, and hence we can find a point
$y \in Y$ such that $x \in \Locus(V;0 \mapsto y)$ by
Theorem~\ref{thm:RC-fibration}. In particular, $x$ and $y$ are
$V$-chain connected. This implies that $y \in G$, and therefore
that $G \cap Y \ne \emptyset$. Since $G$ is a general fiber of
$\f^\o$, this means that $\f^\o(Y \cap X^\o) = S^\o$. Moreover, by the
generality of $G$, we can in fact assume that $G$ has non-empty
intersection with a very general fiber of $\p^\o$. On the other
hand, recalling that $\p$ is defined by the restriction $V_Y$ of
the family $V$ defining $\f$, we see that if $F$ is a very general
fiber of $\p^\o$ meeting $G$, then necessarily $F \subseteq G$. In
conclusion, after possibly further shrinking $Z^\o$ (and
consequently $Y^\o$), we may assume to have a commutative diagram
$$
\xymatrix{
Y^\o \ar[d]_{\p^\o} \ar@{^{(}->}[r] & X^\o \ar[d]^{\f^\o} \\
Z^\o \ar[r]^{\d^\o} & S^\o,
}
$$
where $\d^\o$ is a dominant morphism.

We need to show that, after possibly changing the birational model for
$Y\ratquot_{V_Y}$, the map $\d^\o$ extends to a surjective morphism
$$
\d \colon Y\ratquot_{V_Y} \to X\ratquot_V.
$$
We first take the normalization $Z$ of
a projective compactification of $Z^\o$. Note that $\d^\o$ extends
to a rational map $Z \rat X\ratquot_V$, which is defined by
some linear system $\H$ on $Z$. Then we take
$Z'$ to be the normalization of the blow-up
of the base scheme of $\H$. Since the proper transform of $\H$
on $Z'$ is base-point free, the rational map $\d^\o$
lifts to a well defined morphism $\d \colon Z' \to X\ratquot_V$.
This is a morphism of projective varieties, and $\d^\o$
is dominant, so $\d$ is surjective. Moreover,
the morphism $Z' \to Z$ is an isomorphism over $Z^\o$, and
thus we can identify the latter with a subset of $Z'$.
Then we take $Z'$ as the $V_Y$-quotient $Y\ratquot_{V_Y}$ of $Y$.
This proves the theorem.
\end{proof}

Theorem~\ref{thm:main-RC-extension:part-I} can be viewed
as a general property of rationally connected fibrations.
Aiming also at applications in extension problems of
specific fibrations, that will be addressed in the following sections,
we are interested in determining sufficient conditions to ensure
that the map $\d$, whose existence was proven in the previous
theorem, is generically finite. This is the content of the
next theorem, which is the main result of this section.

\begin{thm}\label{thm:main-RC-extension:part-II}
Let $X$ be a smooth projective variety, and assume that $Y \subset
X$ is a smooth subvariety with ample normal bundle. Let $i \colon
Y \inj X$ be the inclusion, and
suppose that the induced map $i^*\colon N^1(X) \to N^1(Y)$ is surjective.
Let $V = \cup_{\a \in A} V_\a$ be a
family of rational curves on $X$, and assume that there is a
subset $B \subseteq A$ such that
\begin{enumerate}
\item[(i)] $\R_{\ge 0}[V] = \sum_{\b \in B} \R_{\ge 0}[V_\b]$, and
\item [(ii)] $\Locus\big(i_*^{-1}(V_\b)\big)$ is dense in $Y$ for
every $\b \in B$.
\end{enumerate}
Let $V_Y := \;\rangle i_*^{-1}(V) \langle\;$
be the restriction of $V$ to $Y$, and let
$$
\xymatrix{
Y \ar@{-->}[d]_{\p} \ar@{^{(}->}[r] & X \ar@{-->}[d]^{\f} \\
Y\ratquot_{V_Y} \ar[r]^\d & X\ratquot_V,
}
$$
be the commutative diagram given by
Theorem~\ref{thm:main-RC-extension:part-I}.
Suppose that one of the following two conditions holds:
\begin{enumerate}
\item $V$ is an unsplit family; or \item $\codim_XY <
\dim\RC_{V_Y}$.
\end{enumerate}
Then the morphism $\d$ is generically finite.
\end{thm}

\begin{proof}
We first prove that $\d$ is generically finite when the
condition~(a) of the statement is satisfied. If $Y\ratquot_{V_Y}$
is a point, then the statement is obvious, so we can suppose that
$Y\ratquot_{V_Y}$ is positive dimensional.  We suppose by
contradiction that $\d$ is not generically finite.

Keeping the notation introduced in the proof of
Theorem~\ref{thm:main-RC-extension:part-I}, let $G$ be a very
general fiber of $\f^\o$, and consider $Y_G := Y \cap G$. Note that
$G$ is smooth and rationally connected by the restriction $V_G$ of $V$
to $G$ in the sense of Definition \ref{def:restriction-of-V}. So,
since we are assuming that $V$ is unsplit, every irreducible component
of $V_G$ is an unsplit family on $G$, and hence $N_1(G)$ is generated by
classes of $V_G$-curves by Proposition~\ref{prop:unsplit-rho=1}.
Since such curves are all numerically equivalent in $X$,
it follows that the image of the restriction map $\Pic(X) \to \Pic(G)$
has rank 1. Consider the inclusions $j_X
\colon G \inj X$, $j_Y \colon Y_G \inj Y$ and $i_G \colon Y_G \inj
G$. In order to reach a contradiction, we consider the commutative
diagram
$$
\xymatrix{
\Pic(X) \ar[d]_{i^*} \ar[r]^{j_X^*} & \Pic(G) \ar[d]^{i_G^*} \\
\Pic(Y) \ar[r]^{j_Y^*} & \Pic(Y_G).
}
$$
The plan is to give a lower-bound for the rank of the map
$\Pic(X)\to \Pic(Y_G)$ by factoring it through $\Pic(Y)$.

Note that $Y^\o_G := Y^\o \cap G$ is a non-empty
open subset of $Y_G$. Since we are supposing that $\d^\o$ is not generically
finite, we have $\dim\big(\p^\o(Y^\o_G)\big) \ge 1$. Therefore
we can find a curve $\G \subseteq Y_G$ such that
$\G \cap Y^\o_G \ne \emptyset$ and
$\dim\big(\p^\o(\G \cap Y^\o_G)\big) \ge 1$.
Since $\Locus(V_Y)$ is dense in $Y$, every fiber of $\p^\o$ has positive
dimension and, since these fibers are proper,
we can fix a curve $C \subset Y_G$ lying inside a fiber
of $\p^\o|_{Y^\o_G}$.

Now, fix a projective model $Y\ratquot_{V_Y}$ for the rational
quotient of $Y$ containing $Z^\o$ as an open subset, and let $L
\subset Y\ratquot_{V_Y}$ be a general hyperplane section passing
through a point in $\p^\o(\G \cap Y^\o_G)$ but not containing the
point $\p^\o(C)$. Let $L^\o := L|_{Z^\o}$ be the restriction of $L$
to $Z^\o$, and let $D \subset Y$ be the closure (in $Y$) of
$(\p^\o)^{-1}(L^\o)$. If $U \subset Z^\o \setminus L^\o$ is an open
neighborhood of $\p^\o(C)$, then $(\p^\o)^{-1}(U)$ is an open
neighborhood of $C$ in $X^\o$ (hence in $X$) and is disjoint from
$(\p^\o)^{-1}(L^\o)$; this shows that $D \cap C = \emptyset$. Let
then $D_G$ and $H_G$ be the restrictions of $D$ and of a general
hyperplane section $H$ of $Y$ to $Y_G$. By construction, we have
$D_G\.\G > 0$ and $D_G\.C = 0$, whereas $H_G$ is ample; in
particular, $D_G$ and $H_G$ induce numerically independent
elements of $\Pic(Y_G)$. Since both of them are restrictions of
divisors on $Y$, this implies that $\rk \im(j_Y^*) \ge 2$.
Therefore, observing that the cokernel of $i^* \colon \Pic(X) \to
\Pic(Y)$ is torsion due to the surjectivity of $N^1(X) \to
N^1(Y)$, we conclude that
$$
\rk\im(j_Y^*\o i^*)\ge 2.
$$
On the other hand, we have $\rk\im(i_G^*\o j_X^*) \le 1$,
and therefore we have a contradiction by the commutativity of the
above diagram. This proves that $\d$ is generically finite
if the condition~(a) of the statement of the theorem is satisfied.

It remains to prove that $\d$ is generically finite under the
hypothesis~(b). Let $G$ be a very general (smooth and connected)
fiber of $\f^\o$, and let $F$ be a very general fiber of $\p$ among
those that are contained in $G$. Choosing $G$ sufficiently
general, we can also ensure that $F$ is general among the fibers
of $\p^\o$. Then, to conclude the proof of the theorem, we need to
show that
\begin{equation}\label{eq:5}
\dim G = \dim F + \codim_XY.
\end{equation}
Note indeed that,
since $\dim G = \dim \RC_V$ and $\dim F = \dim\RC_{V_Y}$,
$\eqref{eq:5}$ implies that $\d$ is generically finite.

In order to prove \eqref{eq:5}, we fix a very general point $ y
\in F$. By Lemma~\ref{lem:B-dense}, $\Locus(V,\{0\} \mapsto Y)$ is
dense in $X$. Moreover, we know that $\Locus(V,0 \mapsto y)$ is
contained in $G$ by Theorem~\ref{thm:RC-fibration}. Thus
$\Locus(V,0 \mapsto y)$ sweeps a dense subset of $G$, if we let
$y$ vary in $F$. Therefore, by our choices of $F$ and $y$, we can
assume that $\Locus(V,0 \mapsto y)$ contains a very general point
of $G$. Note that any point of $G$ is connected by a $V_G$-chain to the
point $y$, where $V_G$ denotes the restriction of $V$ to $G$ in
the sense of Definition~\ref{def:restriction-of-V}. Therefore
Proposition~\ref{prop:V-chain-connected=Locus(V,y)-dense} applied
to $G$ and $V_G$ implies that $\Locus\big(\|V_G\|, 0 \mapsto
y\big)$ is dense in $G$. Note that the two families $\|V\|(0
\mapsto y)$ and $\|V_G\|(0 \mapsto y)$ are the same since $G$ is
an ${\rm RC}_{\|V\|}$ equivalence class. In particular, we have
$$
\dim G = \dim\Locus\big(\|V\|, 0 \mapsto y\big).
$$

At this point the idea is to replace $V$ with $\|V\|$, so that, by
the above argument, we can directly assume, without loss of
generality, that the two sets $\Locus(V, 0 \mapsto y)$ and $G$
have the same dimension (and, in fact, that the first set is dense
in the second one). In order to make this step, we need to show,
first of all, that $\|V\|$ satisfies hypotheses analogous to (i)
and (ii) imposed to $V$ in the statement of the theorem, and
moreover that replacing $V$ with $\|V\|$ does not affect the
quotient maps $\f$ and $\p$ and the respective rational quotients,
which implies condition (b) for $\|V\|$.

\begin{lem}
The family $\|V\|$ satisfies the hypotheses (i) and (ii) imposed
to $V$ in Theorem~\ref{thm:main-RC-extension:part-II}. Moreover,
the families $V_Y = \;\rangle i_*^{-1}(V)\langle\;$ and $\|V\|_Y
:= \;\rangle i_*^{-1}(\|V\|)\langle\;$ define the same rational
fibration on $Y$. In particular, $\codim_XY < \dim\RC_{\|V\|_Y}$.
\end{lem}

\begin{proof}[Proof of the lemma.]
First recall that $V_Y$ is a numerically covering family by
Theorem \ref{thm:main-RC-extension:part-I}, hence so is $\|V_Y\|$.
We compare $\|V\|$ with
$$
W := \langle i_*(\|V_Y\|_{\rm cov}) \rangle.
$$
Note that $\|V_Y\|_{\rm cov}$ is non-empty and $W$ is a subfamily
of $\|V\|$. We claim that
\begin{equation}\label{eq:4}
\R_{\ge 0}[\|V\|] = \R_{\ge 0}[W] \quad\text{in}\quad N_1(X).
\end{equation}
By the definition of $\|V\|$, this is equivalent to $\R_{\ge 0}[V]
= \R_{\ge 0}[W]$. The inclusion $\R_{\ge 0}[V] \supseteq \R_{\ge
0}[W]$ is obvious, so we need to show that the reverse inclusion
holds. In fact, by the hypothesis (i) on $V$, it suffices to show
that $\R_{\ge 0}[V_\b] \subseteq \R_{\ge 0}[W]$ for every $\b \in
B$. Since $\Locus(i_*^{-1}(V_\b))$ is dense in $Y$, at least one
of the irreducible components of $\;\rangle i_*^{-1}(V_\b)
\langle\;$ is a covering family on $Y$, hence a subfamily of
$\|V_Y\|_{\rm cov}$. This implies that $\R_{\ge 0}[V_\b] \subseteq
\R_{\ge 0}[W]$, hence \eqref{eq:4} is proven.

On the other hand, the restriction to $Y$ of any irreducible
component of $W$ is a union of irreducible components of
$\|V_Y\|_{\rm cov}$, which are covering by Definition
\ref{def:num-covering-family}. Therefore their loci are dense in
$Y$. Combining this with \eqref{eq:4} proves the first assertion.

We are currently assuming that $V$ satisfies the condition~(b) in
Theorem~\ref{thm:main-RC-extension:part-II}. Therefore we get the
desired inequality $\codim_XY < \dim\RC_{\|V\|_Y}$ as soon as we
show that the families $V_Y$ and $\|V\|_Y$ define the same
rational quotient of $Y$, because then $\dim\RC_{V_Y} =
\dim\RC_{\|V\|_Y}$. By using
Proposition~\ref{prop:numerical-condition-on-RC} it is enough to
show that
\begin{equation}\label{eq:1}
\R_{\ge 0}[\|V\|_Y] = \R_{\ge 0}[V_Y] \quad\text{in}\quad N_1(Y).
\end{equation}
The inclusion $\R_{\ge 0} [\|V\|_Y] \supseteq \R_{\ge 0}[V_Y]$ is
obvious. To prove the other, note that the embedding of $Y$ in $X$
induces an inclusion
$$
\iota \colon N_1(Y) \inj N_1(X),
$$
which follows from the hypothesis that $i^*\colon N^1(X) \to
N^1(Y)$ is surjective made in
Theorem~\ref{thm:main-RC-extension:part-II}. Observe that
$$
\iota\Big(\R_{\ge 0} [\|V\|_Y]\Big)
\subseteq \R_{\ge 0}[\|V\|] = \R_{\ge 0}[V] =
\sum_{\b \in B}\R_{\ge 0}[V_\b] = \iota\big(\R_{\ge 0}[V_Y]\big)
\quad\text{in}\quad N_1(X).
$$
Since $\iota$ is injective, the inclusion holds in $N_1(Y)$,
before taking images. This proves equality \eqref{eq:1} and
completes the proof of the lemma.
\end{proof}

We now come back to the proof of the theorem. By the previous lemma,
we are allowed to replace $V$ with $\|V\|$.
Therefore we can directly assume without loss of generality that
\begin{equation}\label{eq:same-dimension-Locus-G}
\dim \Locus(V, 0 \mapsto y) = \dim G.
\end{equation}
Then, in order to show \eqref{eq:5} and hence to conclude the
proof, it suffices to show that $\dim \Locus(V, 0 \mapsto y) =
\dim X - \dim Y\ratquot_{V_Y}$.

For short, let $d = \dim X$ and $k = \dim Y\ratquot_{V_Y}$.
Since $y$ is a general point of $Y$, we can assume that
$V_Y(0 \mapsto y)$ contains a sufficiently general point $[g]$
of a covering component of $V_Y$, so that $g$ is free
and $C := g(\P^1)$ is a smooth curve
(the smoothness of $C$ is not essential for the proof,
but simplifies the notation).

We observe that $N_{Y/X}|_F$ is an ample vector bundle over $F$,
and that $\rk N_{Y/X}|_F < \dim F$ by condition (b), and therefore
$H^1((N_{Y/X}|_F)^*) = 0$ by the Le Potier vanishing theorem
\cite{LeP}. This implies that the short exact sequence
$$
0 \to N_{F/Y} \to N_{F/X} \to N_{Y/X}|_F \to 0
$$
splits since $N_{F/Y} \cong \O_F^{\oplus k}$, and therefore we get a surjection
$N_{F/X} \surj N_{F/Y}$.
Restricting to $C$ and composing with the natural surjections
$T_X|_C \surj N_{C/X}$ and $N_{C/X} \surj N_{F/X}|_C$,
we obtain a chain of surjections
$$
T_X|_C \surj N_{C/X} \surj N_{F/X}|_C \surj N_{F/Y}|_C \cong
\O_{\P^1}^{\oplus k}.
$$
Since $g$ is free, so is $f := i\o g$
by Proposition~\ref{prop:free-on-Y=free-on-X},
and therefore $T_X|_C$ is nef.
Writing $T_X|_C \cong \oplus_i \O_{\P^1}(b_i)$, we have $b_i \ge 0$, and thus
the existence of a quotient of $T_X|_C$ isomorphic to $ \O_{\P^1}^{\oplus k}$
implies that
\begin{equation}\label{eq:bound-on-b_i=0}
\#\{ b_i \mid b_i = 0 \} \ge k.
\end{equation}

Note that $(V,0\mapsto y)$ is smooth at $[f]$ by
Proposition~\ref{prop:free=smooth}. Let
$u\colon\P^1\times (V,0\mapsto y)\to X$
be the universal map, and let $q\in \P^1$ be a point different from $0$.
By \cite[Proposition~II.3.10]{KoBook} and \eqref{eq:bound-on-b_i=0}, we have
$$
\rk du(q,[f]) \le d-k.
$$
Since $(V,0\mapsto y)$ is smooth at $[f]$, we obtain $\dim {\rm
Locus}(V,0\mapsto y) \le d-k$. On the other hand, we have $\dim G
\ge d-k$ by Theorem~\ref{thm:main-RC-extension:part-I}. Thus by
\eqref{eq:same-dimension-Locus-G} we conclude that $\dim {\rm
Locus}(V,0\mapsto y) = d-k$, which proves \eqref{eq:5}. This
concludes the proof of the theorem.
\end{proof}

\begin{rmk}
It is interesting to compare the result of
Theorem~\ref{thm:main-RC-extension:part-II}
when the condition in case~(b) is satisfied
with a result of Sommese \cite[Proposition~III]{Sommese}
(see also \cite[Theorem~(5.3.1)]{Book} for a more general version)
in which the case of an ample divisor endowed with a fibration
is considered. Indeed, one notices that if $\codim_XY = 1$, then
the inequality in~(b) reduces exactly to the hypothesis
imposed in the theorem of Sommese.
\end{rmk}

\begin{rmk}\label{rmk:Hartshorne-conj}
Hartshorne's conjecture \cite[Conjecture 4.5]{Hartshorne} on the
intersection of ample submanifolds
with complementary dimensions, if true, would imply
that the morphism $\d$ in Theorem~\ref{thm:main-RC-extension:part-II}
is in fact an isomorphism whenever the condition given in~(b) is satisfied.
Indeed, let $G$ be a general fiber of $\phi$.
Then $Y \cap G$ is a disjoint union of $\deg\delta$
fibers $F_j$ of $\p$, and by using standard exact sequences one
easily sees that $N_{Y\cap G/G}\cong N_{Y/X}|_{Y\cap G}$. In particular, we have
$N_{F_j/G}\cong N_{Y\cap G/G}|_{F_j}\cong N_{Y/X}|_{F_j}$, which is ample.
Note that condition (b) implies
that $2\dim F_j>\dim G$. Then, assuming that $\deg\d \ge 2$,
Hartshorne's conjecture would give the
contradiction $F_1\cap F_2\neq\emptyset$; therefore $\delta$ can only be
birational. In fact, modulo suitable choices of the rational quotients,
one could even assume that $\delta$ is an isomorphism.
We recall that this conjecture of Hartshorne is known
to be true, for instance, for homogeneous spaces
\cite{Lub}; we will use this fact in the proof of
Theorem~\ref{thm:P-bundle-Q-fibration}.
\end{rmk}

\begin{rmk}\label{rmk:Occhetta-1}
The previous results are related to a result
of Occhetta \cite[Proposition~4]{Occ},
where an analogous extension property is obtained under more
restrictive hypotheses.
In the setting of \cite{Occ}, $Y$ is supposed to be
the zero locus of a regular section of an ample vector bundle on $X$,
$\p$ is the {\rm MRC}-fibration of $Y$, and the rational
quotient of $Y$ is assumed to have positive geometric genus.
\end{rmk}

In applications, it might be useful to start with families of
rational curves on $Y$, rather than only considering restrictions
of families on $X$. The following elementary property will be
useful. Let $X$ and $Y$ be as above,
with the inclusion inducing a surjection $N^1(X) \to N^1(Y)$,
and consider a numerically
covering family $W$ of rational curves on $Y$. Let then $V =
\langle i_*(W) \rangle$ be the family on $X$ obtained as the
extension of $W$, and let $V_Y = \rangle i_*^{-1}(V) \langle$ be
its restriction back to $Y$. Note that $V_Y$ is numerically
covering since $W$ is so. We clearly have $W \subseteq V_Y$, but
in general the inclusion may be strict (one can take as an example
the case of a quadric $Y = \Q^2$ in $X = \P^3$, taking $W$ to be
the family of lines on $Y$ corresponding to one of the two
rulings). By this inclusion and
Theorem~\ref{thm:main-RC-extension:part-I}, we have a commutative
diagram
$$
\xymatrix{
& Y \ar@{-->}[dl]_\t \ar@{-->}[d]^{\p} \ar@{^{(}->}[r] & X \ar@{-->}[d]^{\f} \\
Y\ratquot_W \ar[r]^\g & Y\ratquot_{V_Y} \ar[r]^\d & X\ratquot_V,
}
$$
for a suitable choice of the rational quotients, and some
surjective morphisms $\g$ and $\d$.

\begin{prop}\label{prop:RC-W-versus-V_Y}
Keeping the notation introduced above, both families $W$ and $V_Y$
define the same rational quotient (i.e., $Y\ratquot_W =
Y\ratquot_{V_Y}$ and $\t = \p$) if either one of the following two
situations occurs:
\begin{enumerate}
\item
$\dim Y\ratquot_W = \dim X\ratquot_V$, or
\item
$i^*\colon N^1(X) \to N^1(Y)$ is surjective.
\end{enumerate}
\end{prop}

\begin{proof}
We observe that the general fiber of $\g$ is rationally connected,
since it is dominated by the general fiber of $\p$, and in
particular it is connected. Then the conclusion in case~(a)
follows immediately from the commutativity of the above diagram.
We now focus on case~(b), and thus assume that $N^1(X) \to N^1(Y)$
is surjective. By duality, we have an injection $N_1(Y) \inj
N_1(X)$ and, by the construction of $V_Y$, we deduce that $\R_{\ge
0}[V_Y] = \R_{\ge 0}[W]$. Then the assertion follows from
Proposition~\ref{prop:numerical-condition-on-RC}.
\end{proof}

The case when the rational quotient is one-dimensional is
particularly well suited because of the following well-known property. We are
grateful to the referee for suggesting this proof, which simplifies
our original argument.

\begin{prop}\label{prop:proper-maps-over-curves-are-morphisms}
Let $X$ be a normal projective variety, let
$\f^\o\colon X^\o \to C$ be a dominant morphism
from a non-empty open subset $X^\o\subseteq X$
to a smooth projective curve $C$, and assume that
$\f^\o$ is proper over its image.
Then $\f^\o$ extends to a (surjective) morphism $\f \colon X \to C$.
\end{prop}

\begin{proof}
We consider a projection $C \to \P^1$. Composing with $\phi^\o$,
we have a rational map
$X \rat \P^1$ which is defined by some linear pencil. Because of the properness of
$\phi^\o$, this pencil is free, hence the map $X \to \P^1$ is regular
and we can lift it via Stein factorization
to a regular morphism $X \to C$ which extends $\phi^\o$.
\end{proof}

Coming back to the extension of rationally connected
fibrations, we immediately obtain the following result.

\begin{cor}\label{cor:RC-extension-unsplit-over-curve}
With the same notation and assumptions as in
Theorem~\ref{thm:main-RC-extension:part-II},
suppose furthermore that $\dim Y\ratquot_{V_Y} = 1$.
Then the $\RC_{V_Y}$-fibration $\p$ of $Y$ and the $\RC_V$-fibration $\f$ of $X$
are surjective morphisms fitting in a commutative diagram
$$
\xymatrix{
Y \ar[d]_{\p} \ar@{^{(}->}[r] & X \ar[d]^{\f} \\
Y\ratquot_{V_Y} \ar[r]^\d &X\ratquot_V,
}
$$
where $\d$ is a finite morphism of smooth curves.
\end{cor}

\begin{proof}
It follows from Theorems~\ref{thm:main-RC-extension:part-I}
and~\ref{thm:main-RC-extension:part-II},
and Proposition~\ref{prop:proper-maps-over-curves-are-morphisms}.
\end{proof}

\section{Extension of Mori contractions}

In this section we consider the problem to extend Mori
contractions.
Let $Y$ be a smooth projective variety with Picard number $\r := \r(Y) \ge 2$,
and suppose that $R \subset N_1(Y)$
is a $K_Y$-negative extremal ray of $\CNE(Y)$ such that the
associated Mori contraction
$$
\p \colon Y \to Z
$$
is of fiber type, namely, that $\dim Z < \dim Y$.

We fix a rational curve $C$ in a general fiber of $\p$ passing
through a general point of $Y$, whose numerical class $[C]$
generates $R$. Then we fix a normalization morphism $f\colon \P^1
\to C \subset Y$, and let
$$
W = \langle [f] \rangle \subset \Hom_{\bir}(\P^1,Y)
$$
be the minimal family of rational curves on $Y$ generated by $[f]$.
Note that, by construction,
$$
\R_{\ge 0}[W] = \R_{\ge 0}[C] = R \quad\text{in}\quad N_1(Y),
$$
and $W$ is a covering family on $Y$, since $C$ passes through a
general point of $Y$. Moreover, $\p$ can be considered as the
$\RC_W$-fibration of $Y$.

\begin{thm}\label{thm:extension-Mori-contraction}
With the above notation, assume that $Y$ is a smooth subvariety,
with ample normal bundle, of a smooth projective variety $X$, such
that the inclusion $i \colon Y \inj X$ induces an isomorphism
$i^*\colon N^1(X) \to N^1(Y)$. Let
$$
V := \langle i_*(W) \rangle \subset \Hom_{\bir}(\P^1,X)
$$
be the extension of the family $W$ to $X$, and assume that $V$ is
an unsplit family on $X$.

Then $\R_{\ge 0}[V]$ is a $K_X$-negative extremal ray of $X$.
Moreover, denoting by $\f \colon X \to S$ the Mori contraction of this ray,
there is a finite surjective morphism $\d \colon Z \to S$
giving a commutative diagram
$$
\xymatrix{
Y \ar[d]_{\p} \ar@{^{(}->}[r] & X \ar[d]^{\f} \\
Z \ar[r]^\d & S.
}
$$
\end{thm}

\begin{proof}
By the hypothesis and duality, the inclusion $i \colon Y \inj X$
induces a natural isomorphism $\iota\colon N_1(Y) \cong N_1(X)$.
Let $R_X := \iota(R)$.
This is a ray in $N_1(X)$, and is
contained in $\CNE(X)$, because of the inclusion
$\iota(\CNE(Y)) \subseteq \CNE(X)$.
Moreover, if $C$ is a curve as above, then adjunction formula gives
$$
K_X\.C = (K_Y - \det N_{Y/X})\. C < 0,
$$
and thus $R_X$ is a $K_X$-negative ray.
The main point here is to prove that $R_X$ is an extremal ray of $\CNE(X)$.

We start by observing two things. First of all, we have $R_X =
\R_{\ge 0}[V]$ in $N_1(X)$, and we can consider $\f$ as the
$\RC_V$-fibration of $X$. Moreover, if
$V_Y := \;\rangle i_*^{-1}(V) \langle\;$ is the
restriction of $V$ to $Y$, then we have
$$
\R_{\ge 0}[V_Y] = \iota^{-1}(\R_{\ge 0}[V]) =
\R_{\ge 0}[W] \quad\text{in}\quad N_1(Y),
$$
by the isomorphism $\iota\colon N_1(Y) \cong N_1(X)$. Therefore $V_Y$
and $W$ define the same rational quotient, by
Proposition~\ref{prop:numerical-condition-on-RC}.

Note that $V$ satisfies the conditions~(i) and~(ii)
of Theorem~\ref{thm:main-RC-extension:part-I}.
Then, by Lemma~\ref{lem:B-dense} and the fact that $V$ is unsplit, we have
\begin{equation}\label{eq:2}
\Locus(V,\{0\} \to Y) = X.
\end{equation}
Indeed the lemma implies that $\Locus(V,\{0\} \to Y)$ is dense
in $X$, and the unsplitness of $V$ implies that such locus is proper,
hence closed in $X$.

Let $D$ be any divisor on $X$
whose restriction $D|_Y$ is a good supporting divisor for $R$.
Note that $D\.C = 0$. Let $\G$ be any curve on $X$.
Then, by Proposition~\ref{prop:relative-unspli=picard1:occhetta},
there is a curve $\G_Y \subset Y$ such that
$[\G] =a[\G_Y] + b[C]$ in $N_1(X)$,
with $a \ge 0$. This implies that
$$
D\.\G = a\, D\.\G_Y + b\, D \.C = a\, D|_Y \.\G_Y \ge 0,
$$
since $D\.C = 0$ and $D|_Y$ is nef. This proves that $D$ is nef.

Since $R$ is an extremal ray of $\CNE(Y)$, by duality
it corresponds to it an extremal face of $\Nef(Y)$
of maximal dimension $\r-1$, and therefore we can find $\r-1$ good
supporting divisors of $R$ whose numerical classes are
linearly independent in $N^1(Y)$.
By the previous argument, this implies that such good supporting divisors,
up to numerically equivalence,
extend to divisors on $X$ that are nef and trivial on $R_X$, and whose
numerical classes are linearly independent in $N^1(X)$.
Since there are $\r-1$ of them, and $\r = \r(X)$
by the isomorphism $N^1(X) \cong N^1(Y)$,
this implies that $R_X$ is an extremal ray of $\CNE(X)$.

So, let $\f \colon X \to S$ be the Mori contraction of $R_X$. By taking
the Stein factorization of the restricted morphism
$\f|_Y\colon Y \to S$, we obtain a commutative diagram
$$
\xymatrix{
Y \ar[d]_{\ff} \ar@{^{(}->}[r] & X \ar[d]^{\f} \\
T \ar[r]^\d & S,
}
$$
where $T$ is normal, $\ff$ is surjective with connected fibers,
and $\d$ is finite. We observe that an irreducible curve on $Y$ is
contracted by $\ff$ if and only if it is numerically proportional
(in $X$, and hence in $Y$) to $C$, and therefore if and only if it
is contracted by $\p$. This implies that $T = Z$ and $\ff = \p$.
Moreover, by \eqref{eq:2}, we see that $Y$ meets every fiber of
$\f$, and therefore $\d$ is surjective. This completes the proof
of the theorem.
\end{proof}

\begin{rmk}
Theorem~\ref{thm:extension-Mori-contraction}
can be generalized, with minor changes, to the case of contractions
of $K_Y$-negative extremal faces of $\CNE(Y)$.
\end{rmk}

\begin{rmk}\label{rmk:Occhetta-2}
A result analogous to Theorem~\ref{thm:extension-Mori-contraction}
is obtained in \cite[Proposition 5]{Occ} in a more restrictive
context (in the setting of \cite{Occ}, $Y$ is assumed to be the
zero locus of a regular section of an ample vector bundle on $X$,
with $\dim X\geq 4$, and $\delta$ is proven to be an isomorphism).
Let us note that the proof of \cite[Proposition 5]{Occ} makes use
in an essential way of  the fact that the divisor $D$ (in our and
in his notation as well) is an adjoint divisor, i.e., a divisor of
the form $D=K_X+A$ for some {\em ample} line bundle $A$ on $X$.
This circumstance seems not to be true in general.
The gap was pointed out by Paltin Ionescu. Our arguments fill up
that gap. A similar discussion is given in \cite{BI} in the case $Y$ is
an ample divisor on $X$. In connection with
Theorem~\ref{thm:extension-Mori-contraction} we should also
mention \cite[Theorem~3.4]{AO} and \cite[Lemma~(1.4)]{dFL}.
\end{rmk}

We close this section with the following elementary property on extensions
of relative polarizations, which applies in particular to the
setting of Theorem~\ref{thm:extension-Mori-contraction}.

\begin{prop}\label{prop:extension-of-polarization}
Let
$$
\xymatrix{
Y \ar[d]_\p \ar@{^{(}->}[r]^i & X \ar[d]^\f \\
Z \ar[r] & S
}
$$
be a commutative diagram of morphisms of projective varieties,
and assume that $i$ is an embedding, that $\p$ is not a finite morphism,
and that $\r(X/S) = 1$.
Let $M$ be a line bundle on $X$ whose restriction $M|_Y$ to $Y$ is $\p$-ample.
Then there exists an ample line bundle $H$ on $X$
such that $H|_F \cong M|_F$ for every fiber $F$ of $\p$.
\end{prop}

\begin{proof}
Since the numerical class of a curve in a positive dimensional
fiber of $\p$ spans $N_1(X/S)$, it follows from the
relative Kleiman criterion of ampleness (see, e.g., \cite[Theorem~1.4]{KM})
that $M$ is $\f$-ample. Then the line bundle $H := M + \phi^*(kA)$
is ample if $A$ is an ample line bundle on $S$ and $k \gg 0$
(see, e.g., \cite[Proposition~1.7.10]{Laz}), and $H|_F \cong M|_F$
for every fiber of $\p$.
\end{proof}

\section{Projective bundles and quadric fibrations as ample subvarieties}

The results obtained in the previous sections are here applied
to classify, under suitable conditions, projective manifolds
containing ``ample subvarieties'' that admit a structure of
projective bundles or quadric fibrations.

Here we will adopt the following definitions
of projective bundles and quadric fibrations.
We remark that there are in fact several slightly
different notions of quadric fibrations in the literature,
and the one adopted here is a little more restrictive than others
(see Remark~\ref{rmk:P-bundle-Q-fibrations} below).

\begin{defi}\label{def:P-bundle-Q-fibration}
A surjective morphism $\p\colon Y \to Z$ between smooth projective
varieties is said to be a {\it $\P^m$-bundle} (resp., a {\it
$\Q^m$-fibration}) if $\p$ is an extremal Mori contraction and
there is a line bundle $L$ on $Y$ such that every fiber $F$ of
$\p$ is mapped isomorphically to $\P^m$ (resp., is embedded as an
irreducible and reduced quadric hypersurface of $\P^{m+1}$) via
the complete linear system $|L|_F|$.
\end{defi}

\begin{rmk}\label{rmk:P-bundle-Q-fibrations}
With the above definition, a morphism $\p\colon Y \to Z$ is a
$\P^m$-bundle if and only if $Y \cong \P_Z(\FF)$ for some vector
bundle $\FF$ of rank $m+1$ on $Z$; we can take as $L$ the
tautological line bundle of $\FF$ on $Y$.  On the other hand, not
all scrolls (resp., quadric fibrations) in the adjunction
theoretic sense (cf. \cite[Sections~3.3, 14.1, 14.2]{Book}) are
$\P^m$-bundles (resp., $\Q^m$-fibrations) in our sense. In fact,
our definition of $\Q^m$-fibration does not include conic bundles
with singular fibers, since we require all fibers to be
irreducible; moreover, the assumption that $\p$ is an extremal
Mori contraction implies that $\r(Y/Z) = 1$, which excludes (as we
will see in the proof of
Lemma~\ref{lem:W-irriducible-Q-fibration}) those fibrations with
all fibers isomorphic to $\Q^2$ for which the fundamental group of
the base $Z$ acts with trivial monodromy. Note that the general
fiber of a $\Q^m$-fibration is isomorphic to $\Q^m$.
\end{rmk}

\begin{lem}\label{lem:W-irriducible-Q-fibration}
Let $\p\colon Y \to Z$ be either a $\P^m$-bundle or a $\Q^m$-fibration,
and let $W$ be the family of rational curves on $Y$ generated by the
lines in the fibers of $\p$. Then $W$ is an irreducible family.
\end{lem}

\begin{proof}
The assertion is clear in all cases except when $\p$
is a $\Q^2$-fibration. In this case, a
straightforward extension of the arguments in the proof
of \cite[Proposition~(1.3.1)]{dF2} shows two things. First of all,
$\r(Y/Z) \ne 1$ if and only if all fibers of $\p$ are
isomorphic to $\Q^2$ and the fundamental group of $Z$ acts with
trivial monodromy. Moreover, in all remaining cases, the
lines in the fibers generate an irreducible family.
The fact that the base is one-dimensional in the
setting considered in \cite{dF2} is not essential in the arguments.
\end{proof}

\begin{defi}
An embedding $\P^m \inj \P^n$ is said to be {\it linear} if the image of $\P^m$ is
a linear subspace of $\P^n$. Similarly,
an embedding of irreducible varieties $F \subset G$, with
$F$ isomorphic to a quadric of $\P^{m+1}$ and $G$
isomorphic either to $\P^n$ or to a quadric in $\P^{n+1}$,
is said to be {\it linear} if $\O_G(1)|_F \cong \O_F(1)$.
\end{defi}

For the reminder of this section,
we consider the following setting.
Let $X$ be a smooth projective variety,
and assume that $Y \subset X$ is a positive dimensional
submanifold with ample normal bundle $N_{Y/X}$.
An easy application of a well-known
theorem of Kobayashi and Ochiai gives us the following property,
in which the case when $Y$ is a projective space or a smooth quadric
of dimension $\ge 3$ is considered.

\begin{prop}\label{prop:case-P^m-Q^m}
With the above notation, assume that the inclusion
$i \colon Y \inj X$ induces an isomorphism $i^* \colon \Pic(X) \to \Pic(Y)$,
and denote $n := \dim X$. If $Y \cong \P^m$ for some $m \ge 1$,
then $X \cong \P^n$ and $Y$ is linearly embedded in $X$.
If $Y \cong \Q^m$ for some $m \ge 3$, then either $X
\cong \P^n$ or $X \cong \Q^n$, and $Y$ is linearly embedded in $X$.
\end{prop}

\begin{proof}
Suppose first that $Y \cong \P^m$. By the assumption, we have
$\r(X)=1$, and the ample generator $\O_X(1)$ of $\Pic(X)$
satisfies $\O_X(1)|_Y \cong \O_Y(1)$. Observe that
$$
d := \deg N_{Y/X} \ge \rk N_{Y/X} = n-m.
$$
By adjunction, we have $\O_X(-K_X) \cong \O_X(m+1+d)$,
and therefore $X$ is a Fano manifold of index $\ge n+1$.
Then the assertion follows from Kobayashi and Ochiai's theorem
(see e.g., \cite[(3.1.6)]{Book}).
The argument for the case when $Y \cong \Q^m$ with $m \ge 3$ is analogous.
\end{proof}

\begin{rmk}
The case when $Y \cong \Q^2$ is harder to deal with under such weak
hypotheses. However, if we additionally assume the existence of an
ample line bundle $H$ on $X$ such that $H|_Y \cong
\O_{\P^1\times\P^1}(1,1)$, then this case cannot occur. If at
the same time we relax the hypothesis on the Picard groups, and
only assume that $i^*\colon \Pic(X) \to \Pic(Y)$ is injective,
then it is easy to see that necessarily $\r(X)=1$, either $X \cong
\P^n$ or $X \cong \Q^n$, and $Y$ is linearly embedded in $X$.
Indeed, to exclude the case $\r(X) = 2$, it suffices to apply
Theorem~\ref{thm:extension-Mori-contraction} with the family $W$ being either
one of the two rulings of $\Q^2$. In this way one
obtains two distinct Mori contractions of $X$
whose relative dimensions, added together, exceed the dimension of $X$,
but this is impossible.
\end{rmk}

Let us now recall the following extension of the Lefschetz
hyperplane theorem, essentially due to Sommese
\cite[Proposition~1.16]{Sommese2} (see also
\cite[Theorem~1.1]{LM2}), which will be used in the proof of
case~(b) of Theorem~\ref{thm:P-bundle-Q-fibration} below.

\begin{prop}\label{prop:Lefschetz-Sommese}
Let $X$ be a projective manifold, and suppose that $Y \subset X$
is a submanifold defined by a regular section of an ample vector
bundle on $X$. Suppose that $\dim Y \ge 3$. Then the inclusion $i
\colon Y \inj X$ induces an isomorphism $i^*\colon \Pic(X) \to
\Pic(Y)$.
\end{prop}

We now address the case when $Y$ is a $\P^m$-bundle or a
$\Q^m$-fibration over a smooth projective variety. The following
theorem (which includes Proposition~\ref{prop:case-P^m-Q^m} as a
very special case) is the main result of this section.

\begin{thm}\label{thm:P-bundle-Q-fibration}
Let $Y$ be a submanifold of a projective manifold $X$ with ample normal bundle.
Assume that $Y$ admits either a $\P^m$-bundle or a $\Q^m$-fibration
structure $\p \colon Y \to Z$ over a smooth projective
variety $Z$ (see Definition~\ref{def:P-bundle-Q-fibration}), and that $Y$
has codimension
\begin{equation}\label{eq:hypothesis:m>r}
\codim_XY < \dim Y - \dim Z.
\end{equation}
Further suppose that one of the following two situations occurs:
\begin{enumerate}
\item $Z$ is a curve and the inclusion $i \colon Y \inj X$ induces
an isomorphism $\Pic(X) \cong \Pic(Y)$; or \item $Y$ is the zero
scheme of a regular section of an ample vector bundle $\E$ on $X$.
\end{enumerate}
Then $\p$ extends to a
morphism $\~\p \colon X \to Z$, that is, there is a commutative
diagram
$$
\xymatrix{
Y \ar[d]_\p \ar@{^{(}->}[r]^i & X \ar[dl]^{\~\p} \\
Z.
}
$$
Moreover, letting $n := m + \codim_XY$, the following holds:
\begin{itemize}
\item if $\p$ is a $\P^m$-bundle, then $\~\p$ is a
$\P^n$-bundle; and \item if $\p$ is a $\Q^m$-fibration, then
$\~\p$ is either a $\P^n$-bundle or a $\Q^n$-fibration.
\end{itemize}
In either case, the fibers of $\p$ are linearly embedded in the fibers of $\~\p$.
\end{thm}

\begin{proof}
For clarity of exposition, we discuss case~(a) and case~(b)
separately, as certain steps require different arguments.

\smallskip
Case~(a). Let $W$ be the family of rational curves on $Y$
generated by the lines in the fibers of $\p$. Note that $W$ is an
irreducible family by Lemma~\ref{lem:W-irriducible-Q-fibration},
and that $\p$ is the $\RC_W$-fibration. By
Propositions~\ref{prop:dense=free}
and~\ref{prop:free-on-Y=free-on-X}, $\Hom_{\bir}(\P^1,X)$ is
smooth at the generic point of $i_*(W)$. Let $V = \langle i_*(W)
\rangle$ (note that $V$ is irreducible). By
Proposition~\ref{prop:numerical-condition-on-RC}, the rational
fibration defined by the restriction $V_Y$ of $V$ to $Y$ coincides
with $\p$. Then, by
Corollary~\ref{cor:RC-extension-unsplit-over-curve}, we obtain the
commutative diagram
$$
\xymatrix{
Y \ar[d]_\p \ar@{^{(}->}[r]^i & X \ar[d]^{\f} \\
Z \ar[r]^\d & S,
}
$$
where $\f$ is the $\RC_V$-fibration and $\d$ is a finite morphism
of smooth projective curves.

By definition, there is a line bundle $L$ on $Y$
such that $L|_F \cong \O_F(1)$ for every fiber
$F$ of $\p$. Since $L$ extends to a line bundle on $X$ and $\rho(X/S) = 1$,
Proposition~\ref{prop:extension-of-polarization} implies that
there exists an ample line bundle $H$ on $X$ such that
$H|_F \cong \O_F(1)$ for every $F$.

If $C$ is a line in a fiber of $\pi$ and we set
$b=n$ if $Y$ is a $\P^m$-bundle and $b=n-1$ if $Y$ is a $\Q^m$-fibration,
then we have
$$
(K_X+bH)\.C = (K_Y+bH|_Y)\.C - \det N_{Y/X}\.C < 0.
$$
It thus follows from
the classification of polarized manifolds with large nef-value
due to Fujita and Ionescu
(see \cite[Section~1]{Fujita} or \cite[Section~1]{Ionescu})
that $\phi \colon X \to S$ is either a $\P^n$-bundle or
a fibration with fibers isomorphic to hyperquadrics in $\P^{n+1}$
and relative Picard number 1,
and in fact, since $\dim X \ge 3$ and $\dim S =1$, in the second case
$X$ is a $\Q^n$-fibration.

To conclude, observing that the general fiber of $\phi$ is a homogeneous space,
we deduce from Remark~\ref{rmk:Hartshorne-conj} that $\d$ is in fact
an isomorphism, and therefore we get a diagram
as in the statement by setting $\~\p := \d^{-1}\o\f$.

\smallskip
Case~(b). The first step is to prove the existence of a diagram as
in the statement. We will use different arguments according to the
codimension of $Y$. For short, let $r := \rk\E$. Note that
$\codim_XY = r = n-m$.

If $r = 1$, then \eqref{eq:hypothesis:m>r} implies that $\p$ has
relative dimension $m \ge 2$. Therefore in this case we can apply
a general result of Sommese \cite[Proposition~III]{Sommese}, which
says that $\p$ extends to a morphism $\~\p \colon X \to Z$.

We now assume that $r \ge 2$. Note that in this case $m \ge 3$ by
\eqref{eq:hypothesis:m>r}. To extend $\p$ in this situation, we
will use the results from the previous sections. As in the proof
of case~(a), let $W$ be the irreducible family of rational curves
on $Y$ generated by the lines in the fibers of $\p$, and let $V$
be the unique irreducible component of $\Hom_{\bir}(\P^1,X)$
containing the generic point of $i_*(W)$. By
Proposition~\ref{prop:numerical-condition-on-RC}, the rational
fibration defined by the restriction $V_Y$ of $V$ to $Y$ coincides
with the one defined by $W$, that is, with $\p$. Then, by
Theorems~\ref{thm:main-RC-extension:part-I}
and~\ref{thm:main-RC-extension:part-II}, we obtain a commutative
diagram
$$
\xymatrix{
Y \ar[d]_\p \ar@{^{(}->}[r] & X \ar@{-->}[d]^\f \\
Z \ar@{-->}[r]^\d & S,
}
$$
where $\f$ is the rational quotient defined by $V$ and
$\d$ is a dominant and generically finite map.
Let $G$ be a general fiber of $\f$,
let $F$ be one of the fibers of $\p$ that is contained in $G$,
and fix a line $C \subseteq F$ (we can assume that $G$ and
$F$ are smooth). Note that $\dim G = n$.

\begin{lem}\label{lem:V-unsplit}
Assuming $r \ge 2$, the family $V$ is unsplit.
\end{lem}

\begin{proof}[Proof of the lemma.]
Since the section of $\E$ defining $Y$ in $X$ restricts to a
regular section of $\E|_G$ whose zero scheme is $F$, we can apply
Proposition~\ref{prop:Lefschetz-Sommese} to this setting.
This implies that
the inclusion of $F$ in $G$ induces an isomorphism $\Pic(G) \cong \Pic(F)$.
Note that $\r(F) = 1$, since $F$ is either a projective space or a
quadric of dimension $m \ge 3$. Therefore we have $\r(G) = 1$.
Let $L$ be a line bundle on $Y$ such that $L|_F \cong \O_F(1)$, and let $\~L$
be the extension of $L$ to $X$. We observe that
$\~L|_G$ is ample and $-K_G \equiv a\~L|_G$ for some integer $a$, which is
easily seen to be positive. Then, since $\~L|_G\.C=1$,
we see that $G$ is a Fano manifold of
index $a$. In particular, we have $-K_G\.C = a \le n+1$.
Note that, if $X^\o \subseteq X$ is as in Theorem~\ref{thm:RC-fibration},
then we can assume that $G$ is a fiber of the morphism
$X^\o \to Z$, and hence $K_X|_G = K_{X^\o}|_G = K_G$. Therefore
$$
\det\E\.C = K_Y\.C - K_X\.C = K_F\.C - K_G\.C \le -m+n+1 = r+1.
$$
Since we are assuming that $r \ge 2$, this implies that $\det\E \.C < 2r$.
We conclude that the family $V$ is unsplit by
Remark~\ref{rmk:HC=1-implies-V-unsplit}.
\end{proof}

We come back to the proof of the theorem, still assuming
for the moment that $r \ge 2$. Since $V$ is unsplit
by Lemma~\ref{lem:V-unsplit}, we can apply
Theorem~\ref{thm:extension-Mori-contraction}. This implies that,
in the above diagram, both $\f$ and $\d$ are morphisms and $\d$ is finite,
and moreover $\r(X/S) = 1$.
Note in particular that $\rk \im \big(\Pic(X) \to \Pic(G)\big) = 1$.
We have $\dim G = n$. We set
$b = n$ if $F \cong \P^m$, and $b = n - 1$ if $F \cong \Q^m$.
Denoting by $C$ a line in a general fiber of $\p$, we obtain
$$
(K_G + bH|_G)\cdot C = (K_F + bH|_F)\cdot C - \det N_{F/G}\.C < 0,
$$
since $N_{F/G} \cong N_{Y/X}|_F$ is ample of rank $n-m$. By applying again
\cite{Fujita,Ionescu} and taking into account that $\r(X/S) = 1$, we deduce
that either $(G,H|_G) \cong (\P^n,\O_{\P^n}(1))$ or
$(G,H|_G) \cong (\Q^n,\O_{\Q^n}(1))$, with the second case
occurring only if $F \cong \Q^m$.
Since the general fiber $G$ of $\~\p$ is a homogeneous space,
we deduce from Remark~\ref{rmk:Hartshorne-conj} that $\d$ is in fact
an isomorphism.\footnote{Alternatively, one can deduce this
fact from \cite[Proposition~1.16]{Sommese2}. Indeed this result, applied
to $G\cap Y$ once this is thought of as a subvariety of $G$ defined by a
regular section of $\E|_G$, gives the isomorphism
$H^0(G\cap Y,\Z) \cong H^0(G,\Z)$, and this implies that
$G \cap Y$ is connected.}
We then define $\~\p := \d^{-1}\o\f$.

At this point we have a diagram as in the statement for all values
of $r$. We claim that $\~\p$ is equidimensional with irreducible
and reduced fibers. To see this, let $G_z := \~\p^{-1}(z)$ for an
arbitrary $z\in Z$, and let $\D$ be any component of $G_z$. Note
that $\D\cap Y$ is contained in the fiber $F_z := \p^{-1}(z)$ of
$\p$, which is irreducible and reduced of dimension $m$; in
particular, we have $\dim(\D \cap Y) \le m$. Since $\dim \D \ge n$
and $\D\cap Y$ is defined by a section of $\E|_\D$, we have $\dim
(\D \cap Y) \ge \dim \D - r$. We conclude that $\dim\D=n$ so that
$\~\p$ is equidimensional, and $F_z \subset \D$. Recalling that
$F_z$ is irreducible and reduced, and cut out on $Y$ by $G_z$
(scheme theoretically), we see that $G_z$ is also irreducible and
reduced.
Since the fibers of $\~\p$ are irreducible and reduced, we can
apply the semi-continuity of the $\D$-genus \cite[Section~5
and~(2.1)]{Fujita2}. We conclude that $\~\p$ is either a $\P^n$-bundle or
a $\Q^n$-fibration over $Z$.
\end{proof}

We note that all cases in
Theorem~\ref{thm:P-bundle-Q-fibration} are effective.

\begin{rmk}\label{rmk:Occhetta:3}
We already listed in the introduction
a series of works related to Theorem~\ref{thm:P-bundle-Q-fibration}.
The classical setting when $Y$ is an ample divisor
is widely studied in the literature; for a survey
and complete references we refer to
\cite[Chapter~5]{Book} and \cite{BI}.
In the case when $Y$ has codimension $r \ge 2$,
we are not aware of other results of this type in which the
base of the fibration is allowed to be arbitrary.
\end{rmk}


\begin{thebibliography}{99}

\bibitem{AO}
M. Andreatta and G. Occhetta,
Ample vector bundles with sections vanishing on special varieties.
{\it Internat. J. Math.} {\bf 10} (1999), 677--696.


\bibitem{BI}
M.C. Beltrametti and P. Ionescu,
A view on extending morphisms from ample divisors.
In preparation.

\bibitem{Book}
M.C. Beltrametti and A.J. Sommese,
{\em The Adjunction Theory of Complex Projective Varieties}.
Expositions in Mathematics, 16, W. de Gruyter,  Berlin, (1995).

\bibitem{BSW}
M.C. Beltrametti, A.J. Sommese and J.A. Wi\'sniewski, Results on
varieties with many lines and their applications to adjunction
theory. {\it Complex Algebraic Varieties} (K. Hulek et al. eds.),
Proceedings, Bayreuth, 1990. 16--38, Lecture Notes in Math.,
1507, Springer-Verlag, Berlin, 1992.

\bibitem{Campana}
F. Campana,
Connexit\'e rationnelle des vari\'et\'es de Fano.
{\it Ann. Sci. Ec. Norm. Sup.} {\bf 25} (1992), 539--545.

\bibitem{Debarre}
O. Debarre, {\em Higher-Dimensional Algebraic Geometry}.
Universitext,  Springer-Verlag, Berlin, 2001.

\bibitem{dF1}
T. de Fernex,
Ample vector bundles with sections vanishing along conic fibrations over curves.
{\it Collect. Math.} {\bf 49} (1998), 67--79.

\bibitem{dF2}
T. de Fernex,
Ample vector bundles and intrinsic quadric fibrations over irrational curves.
{\it Matematiche (Catania)} {\bf 55} (2000), 205--222.

\bibitem{dFL}
T. de Fernex and A. Lanteri,
Ample vector bundles and Del Pezzo manifolds.
{\it Kodai Math. J.} {\bf 22} (1999), 83--98.

\bibitem{Fujita2}
T. Fujita,
On the structure of polarized varieties with $\Delta$-genera zero.
{\it J. Fac. Sci. Univ. Tokyo Sect. IA Math.}
{\bf 22} (1975), 103--115.

\bibitem{Fujita}
T. Fujita,
On polarized manifolds whose adjoint bundles are not semipositive.
{\em Algebraic Geometry, Sendai 1985}, Adv. Stud. Pure Math. 10 (1987), 167--178.

\bibitem{GHS}
T. Graber, J. Harris and J.M. Starr,
Families of rationally connected varieties.
{\it J. Amer. Math. Soc.} {\bf 16} (2003), 57--67.

\bibitem{Hartshorne}
R. Hartshorne, {\it Ample Subvarieties of Algebraic Varieties.}
Lecture Notes in Math. 156, Springer-Verlag, Berlin, 1970.

\bibitem{Ionescu}
P. Ionescu,
Generalized adjunction and applications.
{\it Math. Proc. Camb. Philos. Soc.} {\bf 99} (1986), 457--472.

\bibitem{KoBook}
J. Koll\'ar, {\em Rational Curves on Algebraic Varieties}.
Ergeb. Math. Grenzgeb. (3) 32, Springer-Verlag, Berlin, 1996.

\bibitem{KMM}
J. Koll\'ar, Y. Miyaoka and S. Mori,
Rationally connected varieties.
{\it J. Algebraic Geom.} {\bf 1} (1992), 429--448.

\bibitem{KM}
J.~Koll\'{a}r and S.~Mori,
\emph{Birational Geometry of Algebraic Varieties}.
Cambridge Tracts in Mathematics,
Cambridge University Press, Cambridge, 1998.

\bibitem{LM1}
A. Lanteri and H. Maeda,
Ample vector bundles with sections vanishing on projective spaces or quadrics.
{\it Internat. J. Math.} {\bf 6} (1995), 587--600.

\bibitem{LM2}
A. Lanteri and H. Maeda,
Ample vector bundle characterizations of projective bundles
and quadric fibrations over curves.
In Andreatta et al. (eds.),
Proceedings of the international conference: Higher Dimensional
Complex Varieties,
Trento, Italy, June 15--24, 1994. 247--259, W. de Gruyter, Berlin, 1996.

\bibitem{LM3}
A. Lanteri and H. Maeda,
Geometrically ruled surfaces as zero loci of ample vector bundles.
{\it Forum Math.} {\bf 9} (1997), 1--15.

\bibitem{LM4}
A. Lanteri and H. Maeda,
Special varieties in adjunction theory and ample vector bundles.
{\it Math. Proc. Camb. Philos. Soc.} {\bf 130} (2001), 61--75.

\bibitem{Laz}
R. Lazarsfeld, {\it Positivity in Algebraic Geometry. I --
Classical Setting: Line Bundles and Linear Series.} Ergeb. Math.
Grenzgeb. (3) 48, Springer-Verlag, Berlin, 2004.


\bibitem{LeP}
J. Le Potier, Annullation de la cohomologie \`a valeurs dans un
fibr\'e vectoriel holomorphe positif de rang quelconque. {\it
Math. Ann.} {\bf 218} (1975), 35--53.


\bibitem{Lub}
M. L\"ubke,
Beweis einer Vermutung von Hartshorne f\"ur den Fall homogener
Mannigfaltigkeiten.
{\it J. Reine Angew. Math.} {\bf 316} (1980), 215--220.

\bibitem{Occ}
G. Occhetta, Extending rationally connected fibrations.
{\it Forum Math.} {\bf 18} (2006), 853--867.

\bibitem{Sommese}
A.J. Sommese,
On manifolds that cannot be ample divisors.
{\it Math. Ann.} {\bf 221} (1976), 55--72.

\bibitem{Sommese2}
A.J. Sommese,
Submanifolds of Abelian varieties.
{\it Math. Ann.} {\bf 233} (1978), 229--256.

\bibitem{Wis}
J.A. Wi\'sniewski, On a conjecture of Mukai,
{\it Manuscripta Math.} {\bf 68} (1990), 135--141.

\end{thebibliography}
\end{document}